\documentclass[12pt]{article}
\usepackage[dvips]{epsfig}
\usepackage{graphics}
\usepackage{latexsym}
\usepackage{verbatim}
\usepackage{amsmath}
\usepackage{amsthm}
\usepackage{amssymb}

\newtheorem{theorem}{Theorem}[section]
\newtheorem{lemma}[theorem]{Lemma}
\newtheorem{corollary}[theorem]{Corollary}
\newtheorem{proposition}[theorem]{Proposition}

\def\thetheorem {{\arabic{section}.\arabic{theorem}}}

\def\theequation {\arabic{section}.\arabic{equation}}

\def\Remark{\medskip\noindent{\bf Remark: }}
\def\Remarks{\medskip\noindent{\bf Remarks: }}

\def\bul{{$\bullet$\hspace*{2mm}}}

\newcommand{\ens}[1]{\mathbb{#1}}

\newcommand{\N}{\mathbb{N}}

\newcommand{\R}{\mathbb{R}}

\def\cal{\mathcal}
\def\cst{{\rm cst}} %Macro pour constante
\def\RF{{\cal R \cal F}} %Radial Fourier transform
\def\F{{\cal F}} %Fourier transform
\def\Cstab{C_{\rm stab}}
\def\Cdec{C_{\rm dec}}
\def\Creg{C_{\rm reg}}
\def\Cint{C_{\rm int}}
\def\Cbd{C_{\rm BD}}
\def\Cduh{C_{\rm Duh}}

\def\supp{\mbox{supp }}

\def\derpar#1#2{\frac{\partial#1}{\partial#2}}
\def\var{\varepsilon}
\def\signcv{\bigskip\bigskip\hspace{80mm}
\vbox{{\sc C. Villani\par\vspace{3mm}
UMPA, ENS Lyon\par
46 all\'ee d'Italie\par
69364 Lyon Cedex 07\par
FRANCE\par\vspace{3mm}
e-mail:} cvillani@umpa.ens-lyon.fr }}

\def\signcm{\bigskip\bigskip\hspace{80mm}
\vbox{{\sc C. Mouhot\par\vspace{3mm}
UMPA, ENS Lyon\par
46 all\'ee d'Italie\par
69364 Lyon Cedex 07\par
FRANCE\par\vspace{3mm}
e-mail:} cmouhot@umpa.ens-lyon.fr }}

\begin{document}

\title
{Regularity theory for the spatially homogeneous
Boltzmann equation with cut-off}

\author{Cl\'ement Mouhot and C\'edric Villani}

\hyphenation{bounda-ry rea-so-na-ble be-ha-vior pro-per-ties
cha-rac-te-ris-tic parti-cu-lar}

\date{}

\maketitle

\begin{abstract} We develop the regularity theory of the spatially 
homogeneous Boltzmann equation with cut-off and hard potentials
(for instance, hard spheres),
by (i) revisiting the $L^p$-theory to obtain constructive bounds,
(ii) establishing propagation of smoothness and singularities,
(iii) obtaining estimates about the decay of the singularities
of the initial datum. Our proofs are based on a detailed study
of the ``regularity of the gain operator''.
An application to the long-time behavior is presented.
\end{abstract}

\tableofcontents

\section{Introduction}

This paper is devoted to the study of qualitative properties of
solutions to the spatially homogeneous Boltzmann equation with
cut-off and hard potentials. In this work, we shall obtain new,
quantitative bounds on the norms of the solutions in Lebesgue and
Sobolev spaces. Before we explain our results and methods in more
detail, let us introduce the problem in a precise way. 

The spatially homogeneous Boltzmann equation decribes the behavior of
a dilute gas, in which the velocity distribution of particles is
assumed to be independent on the position; it reads
 \begin{equation*}\label{eq:base}
 \derpar{f}{t}  = Q(f,f), \qquad  v \in \R^N, \quad t \geq 0,
 \end{equation*}
where the unknown $f=f(t,v)$ is a time-dependent
probability density on $\R^N$ ($N\geq2$) and $Q$ is the quadratic
Boltzmann collision operator, which we define by the bilinear form
 \begin{equation*}\label{eq:collop}
 Q(g,f) = \int _{\R^N} \, dv_*
 \int _{\ens{S}^{N-1}} \, d\sigma B(|v-v_*|, \cos \theta)(g'_* f' - g_* f).
 \end{equation*}
Here we have used the shorthands $f'=f(v')$, $g_*=g(v_*)$ and
$g'_*=g(v'_*)$, where
 \begin{equation*}\label{eq:rel:vit}
 v' = \frac{v+v_*}{2} + \frac{|v-v_*|}{2} \sigma, \qquad
 v'_* = \frac{v+v^*}{2} - \frac{|v-v_*|}{2} \sigma
 \end{equation*}
stand for the pre-collisional velocities of particles which after
collision have velocities $v$ and $v_*$. Moreover $\theta\in
[0,\pi]$ is the deviation angle between $v'-v'_*$ and $v-v_*$, and
$B$ is the Boltzmann collision kernel (related to the
cross-section $\Sigma(v-v_*,\sigma)$ by the formula $B=\Sigma
|v-v_*|$), determined by physics. On physical grounds, it is
assumed that $B \geq 0$ and that $B$ is a function of $|v-v_*|$ and
$\cos\theta= \left\langle \frac{v-v_*}{|v-v_*|}, \sigma \right \rangle$.

In this paper we shall be concerned with the case when $B$ is {\bf
locally integrable}, an assumption which is usually referred to as
{\em Grad's cut-off assumption} (see~\cite{Grad:58}). The main
case of application is that of hard sphere interaction, where (up
to a normalization constant)
 \begin{equation}\label{eq:collkerhardsphere}
 B (v-v_*, \sigma) = |v-v_*|.
 \end{equation}
We shall study more general kernels than just
\eqref{eq:collkerhardsphere}, but, in order to limit the
complexity of statements, we shall assume that $B$ takes the
simple product form
 \begin{equation}\label{eq:formtenscollker}
 B(v-v_*,\sigma) = \Phi (|v-v_*|) \, b(\cos \theta), \qquad
 \cos \theta = \left\langle \frac{v-v_*}{|v-v_*|}, \sigma \right\rangle.
 \end{equation}
Let us state our assumptions in this context:
 \begin{itemize}
 \item Grad's cut-off assumption, which takes here the simple form
  \begin{equation}\label{eq:hypcut-off}
  \int_0 ^\pi b(\cos \theta) \sin ^{N-2} \theta \, d\theta \ < +\infty.
  \end{equation}
 It is customary in physics and in mathematics to study the case
 when $\Phi$ behaves like a power law $|v-v_*|^\gamma$, and one
 traditionally separates between hard potentials ($\gamma >0$),
 Maxwellian potentials ($\gamma = 0$), and soft potentials ($\gamma
 < 0$). Here we shall concentrate on {\bf hard potentials}, and
 more precisely we shall assume that $\Phi$ behaves like a positive
 power of $|v-v_*|$, in the following sense:
 \item There exists a $\gamma\in (0,2)$ such that
  \begin{equation}\label{eq:hyprad}
  \Phi(0)=0 \ \ \mbox{ and } \ \
  C_\Phi \equiv \left\| \Phi \right\|_{C^{0,\gamma} (\R_+)}  < +\infty.
  \end{equation}
 Here $C^{0,\gamma} (\R_+)$ is the $\gamma$-H\"older space on $\R_+$, i.e
  \begin{equation*}\label{eq:csthold}
  \left\| \Phi \right\|_{C^{0,\gamma} (\R_+)} =
  \sup_{r,s \in \R_+, \ r \neq s} \frac{|\Phi(r) - \Phi(s)|}{|r-s|^\gamma};
  \end{equation*}

 \item In addition to~\eqref{eq:hypcut-off} we shall assume
 a polynomial control on the convergence of the angular integral:
  \begin{equation}\label{eq:hypang2}
  \left|\int_0 ^\pi b(\cos \theta) \sin ^{N-2} \theta \, d\theta
  - \int_\var ^{\pi-\var} b(\cos \theta) \sin ^{N-2} \theta \, d\theta \right|
  \le C_b \var^\delta, \ \ \ \mbox{for some} \ \ \delta >0.
  \end{equation}

\Remark The goal of this assumption is to simplify the computations
and bounds which will be derived. Of course, the $L^1$ integrability of the
angular cross-section implies that the left-hand side in \eqref{eq:hypang2}
goes to~0 as $\var\to 0$, and almost all the results in
the present paper remain true under this sole assumption.
\smallskip \\
 \item Finally, we shall impose a lower bound on the kernel $B$, in
 the form
  \begin{equation} \label{eq:hyplbB}
  \int_{\ens{S}^{N-1}} B(v-v_*,\sigma)\,d\sigma \geq K_B |v-v_*|^\gamma
  \qquad (K_B>0, \quad \gamma>0).
  \end{equation}
 For a kernel in product form, as in~\eqref{eq:formtenscollker}, this
 assumption means that $b$ is not identically (almost everywhere) zero and $\Phi$ satisfies
  \begin{equation} \label{eq:hypradmoins}
  \Phi(|z|) \geq K_\Phi |z|^\gamma \qquad \forall z \in \R^N
  \end{equation}
 for some $K_\Phi>0$.
 \end{itemize}

\Remarks 1. Our assumptions imply that $\Phi$ is bounded from above and below by
constant multiples of $|v-v_*|^\gamma$. In fact, to establish the subsequent
$L^p$ estimates on $Q^+$, it is sufficient to treat this case: since the gain
operator behaves in a monotone way with respect to the collision kernel, the general
estimates follow immediately.
\smallskip

2. It would also be immediate to generalize our results to the case in which
$B$ is a finite sum of products of the form \eqref{eq:formtenscollker}, but much
more tedious to do the same for a general $B$, even if no conceptual difficulty
should arise.

\medskip

The Cauchy problem for hard and Maxwellian potentials is by now
fairly well understood (see for example
Carleman~\cite{Carl:fond:32,Carl:fond:57},
Arkeryd~\cite{Arke:I+II:72}, Mischler and
Wennberg~\cite{MiscWenn:exun:99}, Bobylev~\cite{Boby:maxw:88}),
while soft potentials still remain more mysterious (see
Arkeryd~\cite{Arke:infty:81}, Goudon~\cite{Goud:nocu:97},
Villani~\cite{Vill:nocu:98} for partial results).
\medskip \\
For hard potentials with $0 < \gamma < 2$, the following results are known:
 \begin{itemize}
 \item {\bf Existence and uniqueness} of a solution as soon as the initial
 datum $f_0$ satisfies
  \begin{equation}\label{eq:condinit}
  \int_{\R^N} f_0 (v) \left( 1+|v|^2 \right) \, dv \ < +\infty.
  \end{equation}
 This uniqueness statement in fact holds in the class of solutions with nonincreasing
 kinetic energy, and the solution satisfies the conservation laws
  \begin{equation*}\label{eq:loicons}
  \forall t \geq 0, \ \  \int_{\R^N} f(t,v)
  \left( \begin{array}{c} 1\\v\\ \displaystyle |v|^2 \end{array} \right) \, dv =
  \int_{\R^N} f_0(v)
  \left( \begin{array}{c} 1\\v\\ \displaystyle |v|^2 \end{array} \right) \, dv.
  \end{equation*}
This strong uniqueness result is due to Mischler and
Wennberg~\cite{MiscWenn:exun:99}.
We note that spurious solutions with increasing kinetic energy can be
constructed, see~\cite{Wenn:nonu:99}.

 \item {\bf Boltzmann's $H$-theorem}: let
 $H(f) = \int_{\R^N} f \log f \, dv$, then
 $\frac{d}{dt} H(f(t,\cdot)) \le 0$. In particular, if $H(f_0) < +\infty$, then
  \begin{equation*}\label{eq:Htheo}
  \forall t \geq 0, \qquad H(f(t,\cdot)) \le H(f_0).
  \end{equation*}

 \item {\bf Moment bounds} (Povzner~\cite{Povz:ineq:65},
 Desvillettes~\cite{Desv:momt:93}, Wennberg~\cite{Wenn:momt:94,Wenn:momt:97},
 Mischler and Wennberg~\cite{MiscWenn:exun:99}): if $f_0$ satisfies \eqref{eq:condinit}, then
  \begin{equation*}\label{eq:momentsbounds}
  \forall s \geq 2, \quad \forall t_0 > 0, \qquad \sup _{t \geq t_0} \int_{\R^N}
  f(t,v)(1+|v|^s) \, dv \ < +\infty.
  \end{equation*}
 In words, all moments are bounded for positive times, uniformly as $t$ goes to infinity.
 This effect has been studied at length in the literature, and is strongly linked to the
 behavior of the collision kernel as $|v-v_*| \rightarrow +\infty$. Some explicit bounds
 are available~\cite{Desv:momt:93,Wenn:momt:97}.

\item {\bf Positivity estimates} (Carleman~\cite{Carl:fond:32},
 Pulvirenti and Wennberg~\cite{PulvWenn:binf:97}): without further
 assumptions, it is known that
  \begin{equation*}\label{eq:borneinf}
  \forall t_0 > 0, \ \exists K_0 > 0,\ \exists A_0 > 0; \qquad t\geq t_0 \
  \Longrightarrow \ \forall v \in \R^N, \quad f(t,v) \geq K_0 e^{-A_0 |v|^2}.
  \end{equation*}
 This means that there is immediate appearance of a Maxwellian lower bound (the particles
 immediately fill up the whole velocity space). Again the bounds here are explicit.

\item {\bf $L^p$ bounds}: $L^p$ estimates ($p>1$) have been obtained by several
 authors: Carleman~\cite{Carl:fond:32,Carl:fond:57} and
 Arkeryd~\cite{Arke:Linfty:83} for $p=+\infty$, Gustafsson~\cite{Gust:L^p:86,Gust:L^p:88}
 for $1<p<+\infty$.
 The bounds given by Carleman and Arkeryd are constructive, while this does not
 seem to be the case for Gustafsson's one, obtained by an intricate nonlinear
 interpolation procedure.
 \end{itemize}

Our goal in this work is to complete the picture, while staying in the framework of
hard potentials with cut-off, by
 \begin{itemize}
 \item revisiting the $L^p$ theory ($1<p<+\infty$) and obtain {\bf quantitative
 estimates} together with improved results (holding true
under physically relevant assumptions);

 \item study in detail the phenomena of {\bf propagation of smoothness} and
 {\bf propagation of singularities}, which are certainly the main physical consequences of
 Grad's cut-off assumption.
 \end{itemize}

Unlike Gustafsson's proof, our method does not use the $L^\infty$
theory, nor nonlinear interpolation; it is entirely based on the
important property of ``regularity of the gain operator'', namely
the fact that the positive part of the Boltzmann collision
operator
 \begin{equation*}\label{eq:collopgain}
 Q^+(g,f) = \int _{\R^N}
 \int _{\ens{S}^{N-1}} B\bigl (|v-v_*|, \cos \theta \bigr ) g'_* \, f'\,d\sigma\,dv_*
 \end{equation*}
has a regularizing effect. This phenomenon was discovered by
Lions~\cite{Lion:rgainI+II:94,Lion:rgainIII:94}, and later studied
by Wennberg~\cite{Wenn:rado:94}, Bouchut and
Desvillettes~\cite{BoucDesv:rgain:98}, Lu~\cite{Lu:rgain:98}. On
one hand we shall use some of the results
in~\cite{BoucDesv:rgain:98}, but on the other hand we shall also
need some fine versions of the regularization property which do
not appear in the above-mentioned
references, and this is why we shall devote a whole section to the
study of this regularization effect. This  part should be of
independent interest for researchers in the field, since the $Q^+$
regularity is at the basis of the study of propagation of
regularity for the Boltzmann equation in general, including the
full, spatially inhomogeneous Boltzmann equation.
Wennberg's work~\cite{Wenn:rado:94} will be the starting point of 
our investigation.

Since the pioneering
papers~\cite{Lion:rgainI+II:94,Lion:rgainIII:94} it was known
that the $Q^+$ regularity was useful for smoothness issues; we
shall show here that it is also very powerful for establishing
$L^p$ bounds, as was first suggested in Toscani and
Villani~\cite{ToscVill:cveq:2000}.
In this reference, the case of smoothed soft potentials was
considered; here we shall adapt the strategy to the case of hard
potentials, which will turn out to be much more technical.
Our subsequent study of propagation of smoothness will use these $L^p$ bounds
as a starting point, in the case $p=2$.

Interpolation will play an important role in our estimates, but it will
only be linear interpolation, applied to the bilinear Boltzmann operator 
with one frozen argument (typically, $f\mapsto Q(g,f)$).

Our main results can be summarized as follow: under
assumptions~\eqref{eq:hypcut-off}, \eqref{eq:hyprad},
\eqref{eq:hypang2}, and~\eqref{eq:hyplbB}
 \begin{itemize}
 \item if the initial datum lies in $L^p$, then the solution is bounded in $L^p$, uniformly in time;
 \item if the initial datum is smooth (say in some Sobolev space), then the solution is
 smooth, uniformly in time;
 \item if the initial datum is not smooth, then the solution is not smooth either. However,
 it can be decomposed into the sum of a smooth part (with arbitrary high degree of smoothness)
 and a nonsmooth part whose amplitude decays exponentially fast.
 \end{itemize}
All this will be quantified and stated precisely in
sections~\ref{sec:Lp} and~\ref{sec:propag}. The $L^p$ propagation result is
an improvement of already known results, in the sense that we do not need
extra $L^p$ moment condition on the initial datum; the other results are
new. As an application, we shall
establish some new estimates on the rate of convergence to
thermodynamical equilibrium as time goes to infinity. Although
these estimates are obtained as a consequence of our regularity
study, they will hold true even for nonsmooth solutions.
\medskip

The plan of the present paper is as follows. First, in
section~\ref{sec:prelim}, we give some simple estimates on the
collision operator in various functional spaces. These estimates
will be obtained by simple duality arguments; some of them were
essentially well-known even if maybe not in the particular form
which we give. Then in section~\ref{sec:fine} we begin our fine
study of the regularity of $Q^+$. It is only in
section~\ref{sec:Lp} that we start looking at {\em solutions} of
the Boltzmann equation; in this section we show that if the
initial datum lies in $L^p$ ($1<p<+\infty$) then the solution is
bounded in $L^p$ uniformly in time (besides we prove that a
phenomenon of ``appearance of $L^p$ moments'' occurs, like in the
case $p=1$). In section~\ref{sec:propag}, the main result is a
decomposition theorem of the solution into the sum of a smooth
part (having arbitrary high degree of smoothness) and a nonsmooth
part whose amplitude decays exponentially fast.
As a preliminary we shall also prove propagation of smoothness, and thus
rather precisely tackle the phenomena of propagation of singularities
together with exponential decay. Finally, in
section~\ref{sec:longtime} we give an application to the study of
long-time behavior of the solution: the decomposition theorem
allows one to apply estimates for very smooth solution obtained by
the second author in~\cite{Vill:cveq:TA}, in order to prove rapid
convergence to global equilibrium.

The whole paper is essentially self-contained,
apart from a few simple auxiliary estimates for which precise references
will be given, and from known existence and uniqueness results,
which we here admit. Some facts from linear interpolation theory and harmonic analysis, used
within the proofs, will be recalled in an appendix.

%%%%%%%%%%%%%%%%%%%%%%%%%%%%%%%%%%%%%%%%%%%%%%%%%%%%
%%%%%%%%%%%%%%%%%%%%%%%%%%%%%%%%%%%%%%%%%%%%%%%%%%%%
%%%%%%%%%%%%%%%%%%%%%%%%%%%%%%%%%%%%%%%%%%%%%%%%%%%%

\medskip

\section{Preliminary estimates on the collision operator} \label{sec:prelim}
\setcounter{equation}{0}

Let us first introduce the functional spaces which will be used in
the sequel. Throughout the paper we shall use the notation
$\langle \cdot \rangle = \sqrt{1+|\cdot|^2}$ and we shall denote
by ``$\cst$'' various constants which do not depend on the
collision kernel $B$. Whenever multi-indices are needed we shall
use the common notations $x^\nu = x_1 ^{\nu_1} \cdots x_N
^{\nu_N}$, $\partial^\nu = \partial_1 ^{\nu_1} \cdots \partial_N
^{\nu_N}$, where $\partial_i = \partial/\partial_{x_i}$, and
$\left(
\begin{smallmatrix}  \nu \\ \mu  \end{smallmatrix} \right) =
\left( \begin{smallmatrix}  \nu_1 \\ \mu_1  \end{smallmatrix}
\right)
\cdot \cdot \cdot \left( \begin{smallmatrix}  \nu_N \\
\mu_N  \end{smallmatrix} \right).$
We shall use weighted Lebesgue spaces $L^p_k$ ($p\geq 1$, $k\in\R$)
defined by the norm
 \[ \|f\|_{L^p_k(\R^N)} =  \left ( \int |f(v)|^p
 \langle v \rangle^{pk} dv \right )^{1/p} \] with the convention
 \[ \|f\|_{L^\infty_k(\R^N)} =
 \sup_{v \in \ens{R}^N} \left[ |f(v)| \langle v \rangle^{k}\right].
 \]
We shall also use weighted Sobolev spaces $W^{s,p}_k(\R^N)$;
when $s\in \N$ they are defined by the norm
 \[ \|f\|_{W^{s,p}_k(\R^N)} = \left ( \sum_{|\nu|\leq s}
 \| \partial^\nu f\|_{L^p_k}^p \right )^{1/p}. \]
Then the definition is extended to positive (real) values of $s$ by
interpolation. In particular, we shall denote $W^{s,2}_k = H^s _k$ ;
note that this is a Hilbert space.

We shall make frequent use of the translation operators $\tau_h$ defined by
 \begin{equation*} \label{eq:def:translat}
 \forall \ v \ \in \ \R^N, \ \ \ \tau_h f (v) = f(v-h).
 \end{equation*}
The translation operation does not leave the weighted norms invariant.
Instead, we have the following estimates:
 \begin{equation*} \label{eq:decalage:poids}
 \left\|\tau_h f\right\|_{L^p_{k_1 + k_2}} \le \langle h \rangle^{|k_2|}
 \left\|f\right\|_{L^p_{k_1}}.
 \end{equation*}

Finally, we introduce the $H$ functional:
 \[ H(f) = \int f\log f. \]
For nonnegative functions in $L^1_2$, $H(f)$ is finite if and only if
$f$ belongs to the Orlicz space $L\log L$ defined by the convex
function $\phi(X)= (1+|X|)\log (1+|X|)$.

%%%%%%%%%%%%%%%%%%%%%%%%%%%%%%%%%%%%%%%%%%%%%%%%%%%%
%%%%%%%%%%%%%%%%%%%%%%%%%%%%%%%%%%%%%%%%%%%%%%%%%%%%

\subsection{Some convolution-like inequalities on $Q^+$}

In this subsection we prove some estimates on $Q^+$ in Lebesgue and
Sobolev spaces. In the case of Lebesgue spaces, they are essentially
contained in Gustafsson~\cite{Gust:L^p:86,Gust:L^p:88}; but our
method, based on duality, provides somewhat simpler proofs.

We shall establish two different types of estimates: for the bilinear
Boltzmann collision operator on one hand, and for the quadratic operator
on the other hand. To establish the bilinear estimates, we shall
impose an additional assumption on the angular kernel:
{\em no frontal collision should occur}, i.e. $b(\cos\theta)$
should vanish for $\theta$ close to $\pi$:
 \begin{equation}\label{eq:hypangconvol}
 \qquad \exists \, \theta_b>0 \, ; \quad
 \supp \, b \, ( \cos \theta) \subset \left\{ \theta \ /
 \ \ 0 \le \theta \le \pi - \theta_b \right\}
 \end{equation}
This additional assumption will not be needed, on the other hand,
for the quadratic estimates, i.e. the estimates on $Q^+(f,f)$. Indeed,
$Q^+(f,g) = \tilde{Q}^+(g,f)$ if $\tilde{Q}^+$ is a Boltzmann
gain operator associated with the kernel $\tilde{b}(\cos \theta)
= b(\cos (\pi - \theta))$. In particular, $b(\cos\theta)$ and
$[b(cos\theta) + b(\cos(\pi-\theta))]1_{\cos\theta\geq 0}$ define
the same quadratic operator $Q^+$, and the latter
satisfies~\eqref{eq:hypangconvol} automatically. We note that
$Q^+(g,f)$ and $Q^+(f,g)$ will not necessarily satisfy the same
estimates, since assumption~\eqref{eq:hypangconvol} is not symmetric.
To exchange the roles of $f$ and $g$, we will therefore be led
to introduce the assumption that no grazing collision should
occur, i.e.
 \begin{equation}\label{eq:hypangconvol2}
 \qquad \exists \, \theta_b>0 \, ; \quad
 \supp \, b \, ( \cos \theta) \subset \left\{ \theta \ /
 \ \ \theta_b \le \theta \le \pi \right\}.
 \end{equation}

 \begin{theorem}\label{theo:conv}
 Let $k,\eta \in \R$, $s\in\R_+$, $p\in[1, +\infty]$, and let $B$ be a
 collision kernel of the form~\eqref{eq:formtenscollker},
 satisfying the assumption~\eqref{eq:hypangconvol}.
 Then, we have the estimates
   \begin{equation}\label{eq:convol}
   \left\|Q^+(g,f)\right\|_{L^p _{\eta} (\R^N)} \le C_{k,\eta,p} (B)
   \left\|g\right\|_{L^1 _{|k+\eta|+|\eta|} (\R^N)}
   \left\|f\right\|_{L^p _{k+\eta} (\R^N)},
   \end{equation}
   \begin{equation}\label{eq:dual:fort:deriv}
   \left\|Q^+(g,f)\right\|_{W^{s,p} _{\eta} (\R^N)} \le C_{k,\eta,p} (B)
   \left\|g\right\|_{W^{\lceil s \rceil,1} _{|k+\eta|+|\eta|} (\R^N)}
   \left\|f\right\|_{W^{s,p} _{k+\eta} (\R^N)},
   \end{equation}
 where $C_{k,\eta,p} (B) = \cst \
 (\sin (\theta_b /2))^{\min(\eta, 0)-2/p'}
 \left\|b\right\|_{L^1(\ens{S}^{N-1})}
 \left\|\Phi\right\|_{L^\infty _{-k}}$.
 If on the other hand assumption~\eqref{eq:hypangconvol} is
 replaced by assumption~\eqref{eq:hypangconvol2}, then the same
 estimates hold with $Q^+(g,f)$ replaced by $Q^+(f,g)$.
 \end{theorem}

 \begin{corollary} \label{coro:conv}
 Let $k,\eta\in\R$, $p\in [1,+\infty]$,
 and let $B$ be a collision kernel of the form~\eqref{eq:formtenscollker}.
 Then we have the estimates
  \begin{equation*}\label{eq:convolsym}
  \left\|Q^+(f,f)\right\|_{L^p _{\eta} (\ens{R}^N)} \le C_{k} (B)
  \left\|f\right\|_{L^1 _{|k+\eta|+|\eta|} (\ens{R}^N)}
  \left\|f\right\|_{L^p _{k+\eta} (\ens{R}^N)},
  \end{equation*}
  \begin{equation}\label{eq:dual:sym:deriv}
  \left\|Q^+(f,f)\right\|_{W^{s,p} _{\eta} (\ens{R}^N)} \le C_{k} (B)
  \left\|f\right\|_{W^{\lceil s \rceil,1} _{|k+\eta|+|\eta|} (\ens{R}^N)}
  \left\|f\right\|_{W^{s,p} _{k+\eta} (\ens{R}^N)},
  \end{equation}
 where $C_{k} (B) = \cst \ \left\|b\right\|_{L^1(\ens{S}^{N-1})}
 \left\|\Phi\right\|_{L^\infty _{-k}} $.
 \end{corollary}

\Remarks 1. Of course, if $B$ satisfies assumption~\eqref{eq:hyprad},
then $C_{k}(B)$ is finite as soon as $k\geq \gamma$.
\smallskip

2. No regularity is needed on the collision kernel here.
\smallskip

3. In the particular case $\eta \geq 0$, it is possible
to obtain slightly better weight exponents in Theorem~\ref{theo:conv}
and Corollary~\ref{coro:conv}. One can indeed use the inequality
 \[ |v|^2 \le |v'|^2 + |v' _*|^2 \]
to split the weight on the two arguments of $Q^+$ and get
 \[ \left\|Q^+(g,f)\right\|_{L^p _{\eta}} \le \cst
 \left\|Q^+(G,F)\right\|_{L^p} \]
where $F(v) = f(v) \langle v \rangle^\eta$ and
$G(v) = g(v) \langle v \rangle^\eta$. When $\eta \geq 0$,
the conclusion of Theorem~\ref{theo:conv} thus becomes
 \begin{equation*} \label{eq:convol:poids>=0}
 \left\|Q^+(g,f)\right\|_{L^p _{\eta} (\R^N)} \le C_{k,\eta,p} (B)
 \left\|g\right\|_{L^1 _{k+\eta} (\R^N)}
 \left\|f\right\|_{L^p _{k+\eta} (\R^N)},
 \end{equation*}
and
 \begin{equation*} \label{eq:dual:fort:deriv:poids>=0}
 \left\|Q^+(g,f)\right\|_{W^{s,p} _{\eta} (\R^N)} \le C_{k,\eta,p} (B)
 \left\|g\right\|_{W^{\lceil s \rceil,1} _{k+\eta} (\R^N)}
 \left\|f\right\|_{W^{s,p} _{k+\eta} (\R^N)}
 \end{equation*}
\smallskip

4. As we said above, the corollary is obtained from the theorem upon
replacing $b(\cos\theta)$ by $[b(\cos\theta)+b(-\cos\theta)]
1_{0\leq\theta\leq \pi/2}$. We note that in the case of hard-sphere
collision kernel, the physically relevant regime is $\cos\theta\leq 0$,
so our trick to reduce to $\cos\theta\geq 0$ should just be considered
as a mathematical convenience (which could have been avoided by chosing
different conventions; however there is some other motivation for our
present conventions).
%NdCV: je pense surtout a la coherence avec le cas sans cut-off, comme on l'a
%developpe avec Desvillettes, et dans mon survey.
\medskip

\begin{proof}[Proof of Theorem \ref{theo:conv}]
By duality,
 \begin{equation*} \label{eq:ecriture:norme:duale}
 \left\| Q^+(g,f)\right\|_{L^p _\eta} =
 \sup \left\{ \int Q^+ (g,f) \psi \ \ ; \ \left\|\psi\right\|_{L^{p'} _{-\eta}} \le 1 \right\}.
 \end{equation*}
We apply the well-known pre-postcollisional change of variables,
namely $(v,v_*,\sigma) \to (v',v'_*, (v-v_*)/|v-v_*|)$, which has
unit Jacobian, to obtain
 \begin{equation*} \label{eq:dual}
 \int_{\R^N} Q^+ (g,f) \psi \, dv = \int_{\R^{2N}} dv \, dv_* \, g_* f
 \left( \int_{\ens{S}^{N-1}} B(|v-v_*|, \sigma) \psi (v') \, d\sigma \right)
 \end{equation*}
for all $\left\|\psi\right\|_{L^{p'} _{-\eta}} \le 1$. Let us
define the linear operator $S$ by
 \begin{equation*} \label{eq:opS}
 S\psi(v) = \int_{\ens{S}^{N-1}} B(|v|,\sigma) \psi
 \left( \frac{v + |v| \sigma}{2} \right) d\sigma .
 \end{equation*}
Then
 \begin{equation} \label{eq:lien:Q+:S}
 \int_{\R^N} Q^+ (g,f) \psi \, dv =
 \int_{\R^N} g(v_*) \left( \int_{\R^N} f(v)
 \left(\tau_{v_*} S(\tau_{-v_*} \psi)\right) (v)\,dv \right)\,dv_*.
 \end{equation}
We shall study the operator $S$ in weighted $L^1$ and $L^\infty$ norms.
For brevity we denote
$v^+ = \left( \frac{v + |v| \sigma}{2} \right)$.
By use of the inequality
 \begin{equation*}
 \sin \left(\frac{\theta_b}{2}\right) \, |v| \le |v^+| \le |v|
 \end{equation*}
which is a consequence of~\eqref{eq:hypangconvol}, we find
 \begin{equation}\label{eq:ineq:L^infty}
 \left\|S\psi\right\|_{L^\infty _{-k-\eta}} \le \cst \ (\sin
 (\theta_b /2))^{\min(\eta, 0)} \left\|b\right\|_{L^1(\ens{S}^{N-1})}
 \left\|\Phi\right\|_{L^\infty _{-k}} \left\|\psi\right\|_{L^\infty _{-\eta}}.
 \end{equation}

Next, we turn to the $L^1$ estimate. First,
 \begin{multline*}
 \left\|S\psi\right\|_{L^1 _{-k -\eta}} = \int_{\R^N}
 \int_{\ens{S}^{N-1}} \Phi(|v|) \langle v \rangle ^{-k -\eta} \, b(\cos \theta)
 |\psi(v^+)| \, d\sigma \, dv \\
 \leq (\sin (\theta_b/2))^{\min(\eta, 0)}
 \left\|\Phi\right\|_{L^\infty _{-k}}
 \int_{\R^N} \int_{\ens{S}^{N-1}} \, b(\cos \theta) \,
 |\psi(v^+)| \langle v^+ \rangle ^{-\eta} \, d\sigma \, dv
 \end{multline*}
The change of variable $v \rightarrow v^+$ is allowed
because $b$ has compact support in $[0,\pi-\theta_b]$, and its
Jacobian is $\frac{2^{N-1}}{\cos^2\theta/2}$. By applying it we find
 \begin{eqnarray}\label{eq:ineq:L^1}
 \left\|S\psi\right\|_{L^1 _{-k -\eta}}
 &\le& \cst \, (\sin (\theta_b/2))^{\min(\eta, 0)}
 \left\|\Phi\right\|_{L^\infty _{-k}}
       \int_{\R^N} \int_{\ens{S}^{N-1}} \, b(\cos \theta) \,
       |\psi(v^+)| \langle v^+ \rangle ^{-\eta} \, d\sigma \, dv \\ \nonumber 
 &\le& \cst (\sin (\theta_b/2))^{\min(\eta, 0)} \left\|\Phi\right\|_{L^\infty _{-k}}
       \int_{\R^N} \int_{\ens{S}^{N-1}} b(\cos \theta)
       |\psi(v^+)| \langle v^+ \rangle^{-\eta}
       \frac{2^{N-1}}{\cos^2\theta/2} \, dv^+ \, d\sigma \\ \nonumber
 &\le& \cst \ (\sin (\theta_b/2))^{\min(\eta, 0)-2}
 \left\|b\right\|_{L^1(\ens{S}^{N-1})} \left\|\Phi\right\|_{L^\infty _{-k}} \
 \left\|\psi\right\|_{L^1 _{-\eta}}.
 \end{eqnarray}
By the Riesz-Thorin interpolation theorem (see Appendix), from
inequalities~\eqref{eq:ineq:L^infty} and~\eqref{eq:ineq:L^1} we deduce
 \begin{equation*} \label{eq:ineq:L^p}
 \left\|S\psi\right\|_{L^p _{-k-\eta}} \le C_{k,\eta,p'} (B) \
 \left\|\psi\right\|_{L^p _{-\eta}}, \ \ 1 \le p \le \infty
 \end{equation*}
where $C_{k,\eta,p'} (B) = \cst \ (\sin ( \theta_b/2))^{\min(\eta, 0)-2/p}
\left\|b\right\|_{L^1(\ens{S}^{N-1})} \left\|\Phi\right\|_{L^\infty _{-k}}$.
Plugging this inequality in~\eqref{eq:lien:Q+:S}, we find
 \begin{eqnarray*}
 \left| \int_{\R^N} Q^+ (g,f)\, \psi\, dv \right| &\le&
 \int_{\R^N} dv_* |g_*| \left( \int_{\R^N} dv |f|
 |\left( \tau_{-v_*} S(\tau_{v_*} \psi) \right) (v)|
 \right) \\ \nonumber
 &\le& \int_{\R^N} dv_* |g_*| \left\|f\right\|_{L^p _{k+\eta}}
 \left\|\tau_{-v_*} S(\tau_{v_*} \psi)\right\|_{L^{p'} _{-k-\eta}} \\ 
 &\le& \left\|f\right\|_{L^p _{k+\eta}} \int_{\R^N} |g_*| \langle v_* \rangle^{|k+\eta|}
 \left\|S(\tau_{v_*} \psi)\right\|_{L^{p'} _{-k-\eta}} \, dv_* \\ \nonumber
 &\le& C_{k,\eta,p} (B) \left\|f\right\|_{L^p _{k+\eta}}
 \int_{\R^N} |g_*| \langle v_* \rangle^{|k+\eta|}
 \left\|\tau_{v_*} \psi\right\|_{L^{p'} _{-\eta}} \, dv_* \\ \nonumber
 &\le& C_{k,\eta,p} (B) \left\|f\right\|_{L^p _{k+\eta}} \left\|\psi\right\|_{L^{p'} _{-\eta}}
 \int_{\R^N} |g_*| \langle v_* \rangle^{|k+\eta|+|\eta|} \, dv_*  \\ \nonumber
 &\le& C_{k,\eta,p} (B) \left\|f\right\|_{L^p _{k+\eta}}
 \int_{\R^N} |g_*| \langle v_* \rangle^{|k+\eta|+|\eta|} \, dv_* \\ \nonumber
 &\le& C_{k,\eta,p} (B) \left\|f\right\|_{L^p _{k+\eta}} \left\|g\right\|_{L^1 _{|k+\eta|+|\eta|}}
 \end{eqnarray*}
This concludes the proof of~\eqref{eq:convol}.

We now turn to the proof of~\eqref{eq:dual:fort:deriv}.
It is based on the formula
 \begin{equation}  \label{lem:der:bil}
 \nabla Q^\pm(g,f) = Q^\pm (\nabla g,f) + Q^\pm(g, \nabla f)
 \end{equation}
which is an easy consequence of the bilinearity and the Galilean 
invariance property of the Boltzmann operator, namely
$\tau_h Q(g,f)=Q(\tau_h g, \tau_h f)$. From~\eqref{lem:der:bil}
one can easily deduce a Leibniz formula for derivatives of $Q^+$
at any order, and equation~\eqref{eq:dual:fort:deriv} easily
follows for any $s\in \N$. Indeed, whenever $s \in \N$ we can
apply Theorem \ref{theo:conv} to each term of the Leibniz formula
for $\partial^\nu Q^+ (g,f)$ and find
 \begin{eqnarray*}
 \left\|Q^+(g,f)\right\|^p _{W^{s,p} _{\eta}} &=&
 \sum _{|\nu| \le s} \left\|\partial^\nu Q^+(g,f)\right\|^p _{L^p _{\eta}} \\ 
 &=& \sum _{|\nu| \le s} \ \sum_{\mu \le \nu} \begin{pmatrix}
 \nu \\ \mu \end{pmatrix}
 \left\| Q^+(\partial^\mu g,\partial^{\nu-\mu} f)\right\|^p _{L^p _{\eta}} \\ 
 &\le& C_{k,\eta,p} (B) \sum _{|\nu| \le s} \
 \sum_{\mu \le \nu} \begin{pmatrix}
 \nu \\ \mu \end{pmatrix}
 \left\|\partial^\mu g\right\|^p _{L^1 _{|k+\eta|+|\eta|}}
 \left\|\partial^{\nu-\mu} f\right\|^p _{L^p _{k+\eta}} \\ 
 &\le& C_{k,\eta,p} (B) \left\|g\right\|^p _{W^{s,1} _{|k+\eta|+|\eta|}}
 \left\|f\right\|^p _{W^{s,p} _{k+\eta}}.
 \end{eqnarray*}
Then the general case of~\eqref{eq:dual:fort:deriv} is obtained
by use of the Riesz-Thorin interpolation theorem, with respect to the
variable $f$.
\end{proof}

%%%%%%%%%%%%%%%%%%%%%%%%%%%%%%%%%%%%%%%%%%%%%%%%%%%%
%%%%%%%%%%%%%%%%%%%%%%%%%%%%%%%%%%%%%%%%%%%%%%%%%%%%

 \subsection{A lower bound on $Q^-$}

We shall use the following estimates on $Q^-$.
 \begin{proposition}\label{prop:minor}
 Assume that the collision kernel $B$ satisfies~\eqref{eq:hyplbB}.
 Then, for all $f \in L^1_2$ with $H(f)<+\infty$, there exists a
 constant $K(f)$, only depending on a lower bound on $\int f\,dv$,
 and upper bounds on $\int f|v|^2\,dv$ and $H(f)$, such that
  \begin{equation}\label{eq:lowbdQ^-}
  Q^-(f,f) \geq K(f) \> f(v)\, (1 + |v|)^\gamma.
  \end{equation}
 Similarly, if $|v-v_*|^\gamma$ in the right-hand side
 of~\eqref{eq:hyplbB} is replaced by $\min(|v-v_*|^\gamma,1)$, then
 the conclusion~\eqref{eq:lowbdQ^-} should be replaced by
  \begin{equation}\label{eq:lowbdQ^-:faible}
  Q^-(f,f) \geq K(f) \> f(v).
  \end{equation}
 \end{proposition}

The result is well-known: see for
instance~\cite[lemma~4]{Arke:Linfty:83},
or~\cite[lemma~6]{DesvVill:LChomII:2000}.

%%%%%%%%%%%%%%%%%%%%%%%%%%%%%%%%%%%%%%%%%%%%%%%%%%%%
%%%%%%%%%%%%%%%%%%%%%%%%%%%%%%%%%%%%%%%%%%%%%%%%%%%%
%%%%%%%%%%%%%%%%%%%%%%%%%%%%%%%%%%%%%%%%%%%%%%%%%%%%

\section{Regularity of the gain operator} \label{sec:fine}
\setcounter{equation}{0}

It is known since the works of
P.-L.~Lions~\cite{Lion:rgainI+II:94} that, under adequate
assumptions on the collision kernel $B$, the gain operator
$Q^+(g,f)$ acts like a regularizing operator on each of its
components when the other one is frozen. In this section we shall
establish various versions of this regularizing effect. The
results will of course depend on the assumptions imposed on $B$.

The proof in \cite{Lion:rgainI+II:94} was very technical;
it relied on Fourier integral operators, and the theory of
generalized Radon transform (integration over a moving
family of hypersurfaces), which was studied in detail by Sogge and
Stein at the end of the
eighties~\cite{SoggStei:radon:85,SoggStei:radon:86,SoggStei:radon:90}.
Later Wennberg~\cite{Wenn:rado:94} simplified the proof by using the
Carleman representation~\cite{Carl:fond:32} of $Q^+$, and classical
Fourier transform tools. Both authors prove functional inequalities
which are roughly speaking of the type
 \begin{equation} \label{eq:type1}
 \|Q^+(g,f)\|_{H^{(N-1)/2}} \leq C \|f\|_{L^2}\|g\|_{L^1}.
 \end{equation}
A slightly different family of inequalities was obtained by
much simpler means in independent papers by Bouchut and
Desvillettes~\cite{BoucDesv:rgain:98} and Lu~\cite{Lu:rgain:98}: they
established functional inequalities of the type
 \begin{equation}\label{eq:type2}
 \left\|Q^+ (f,f)\right\|_{H^{(N-1)/2}} \le
 C \left \|f\right\|_{L^2}^2.
 \end{equation}
Four our purposes in the next section, inequalities of
type~\eqref{eq:type2} will not be sufficient, and we shall need
the full strength of inequalities of type~\eqref{eq:type1}. On the
other hand, formulas of the type of~\eqref{eq:type2} will be
sufficient for our regularity study later in the paper.

The precise variants of~\eqref{eq:type1} which will be used in the
sequel cannot be found in~\cite{Wenn:rado:94}, so we shall
re-establish them from scratch. Our proof follows essentially the
idea of Wennberg~\cite{Wenn:rado:94}, and our main contributions
will be to make the constants depend more explicitly on the
features of the collision kernel, to extend the results to
weighted Sobolev spaces of arbitrary order and arbitrary weight,
and to extend the range of admissible collision kernels, allowing
a possible deterioration of the exponents of regularization. It
would also be possible to adapt the proofs by Sogge and Stein,
which are more systematic; but it would be much more tedious to
keep track of the constants.

%Je supprime la rq suivante qui me parait obscure:
%
%There is an important difference to note between inequalities~\eqref{eq:type1}
%and~\eqref{eq:type2}: in the first case,
%thanks to the strong smoothness assumption on the
%kernel, the gain operator behaves like a bilinear
%function of its collision kernel and one of its components
%according to the differentiation.

%%%%%%%%%%%%%%%%%%%%%%%%%%%%%%%%%%%%%%%%%%%%%%%%%%%%
%%%%%%%%%%%%%%%%%%%%%%%%%%%%%%%%%%%%%%%%%%%%%%%%%%%%

 \subsection{A splitting of $Q^+$}

We shall first prove the regularity property on
the gain operator when the collision kernel is very smooth. 
Then we shall include the non-smooth part of the kernel,
at the price of deteriorating the exponents, by an interpolation
procedure with the convolution-like
inequalities of section~\ref{sec:prelim}. This interpolation is
not needed for the proof of propagation of the $L^p$-bound but
will be useful for the study of the propagation of
singularity/regularity performed in section~\ref{sec:propag}.
This calls for an appropriate splitting of the collision kernel,
and therefore of the gain operator.

Let us consider a collision kernel $B= \Phi \, b$ satisfying the
general assumptions~\eqref{eq:hypcut-off},~\eqref{eq:hyprad},
~\eqref{eq:hypang2},~\eqref{eq:hyplbB}. Let
$\Theta: \R \rightarrow \R_+$ be an even $C^\infty$ function such
that $\supp \Theta \subset (-1,1)$, and $\int_\R \Theta \, dx = 1$
and $\widetilde{\Theta}: \R^N \rightarrow \R_+$ be a radial
$C^\infty$ function such that $\supp \widetilde{\Theta} \subset
B(0,1)$ and $\int_{\R^N} \widetilde{\Theta} \, dx = 1$. Introduce
the regularizing sequences
 \begin{equation*}
 \left\{
  \begin{array}{l}
  \Theta_m (x) = m \, \Theta(mx) \qquad (x\in\R) \\ \\
  \widetilde{\Theta}_n (x) =
  n^N \widetilde{\Theta}(n x) \qquad (x\in\R^N).
  \end{array}
 \right.
 \end{equation*}
We shall use these mollifiers to split the collision kernel into a
smooth and a non-smooth part. As a convention, we shall use
subscripts $S$ for ``smooth'' and $R$ for ``remainder''. First, we
set
 \[\Phi_{S,n} = \widetilde{\Theta}_n \ast
 \left( \Phi \ 1_{\ens{A}_n} \right), \qquad
 \Phi_{R,n} = \Phi - \Phi_{S,n}, \]
where $\ens{A}_n$ stands for the annulus $\ens{A}_n =
\left\{ x \in \R^N \ ; \ \frac{2}{n} \le |x| \le n \right\}$.
Similarly, we set
 \[ b_{S,m} = \Theta_m \ast \left( b \ 1_{\ens{I}_m} \right), \qquad
 b_{R,m} = b - b_{S,m},\]
where $\ens{I}_m$ stands for the interval $\ens{I}_m = \left\{ x
\in \R \ ; \ -1+\frac{2}{m} \le |x| \le 1-\frac{2}{m} \right\}$
($b$ is understood as a function defined on $\R$ with compact
support in $[-1,1]$). Finally, we set
 \begin{equation*} \label{eq:splitQ^+}
 Q^+ = Q^+ _{S} + Q^+ _{R}
 \end{equation*}
where
 \begin{equation}\label{eq:splitQ^+S}
 Q^+ _{S}(g,f) = \int _{\R^N} \, dv_*
 \int _{\ens{S}^{N-1}} \, d\sigma \, \Phi_{S,n}(|v-v_*|) \,
 b_{S,m} (\cos \theta) \, g'_* \, f'
 \end{equation}
and
 \begin{equation*} \label{eq:splitQ^+R}
 Q^+ _{R} = Q^+ _{RS} + Q^+ _{SR} + Q^+ _{RR}
 \end{equation*}
with the obvious notations
 \begin{equation*} \label{eq:splitQ^+R:detail}
 \begin{cases} \displaystyle\begin{array}{rcl}
Q^+ _{RS} (g,f) &=& \displaystyle\int _{\R^N} \, dv_*
 \int _{\ens{S}^{N-1}} \, d\sigma \, \Phi_{R,n} \, b_{S,m} \, g'_* \, f'\\ \\
\displaystyle Q^+ _{SR} (g,f) &=& \displaystyle\int _{\R^N} \, dv_*
 \int _{\ens{S}^{N-1}} \, d\sigma \, \Phi_{S,n} \, b_{R,m} \, g'_* \, f'  \\ \\
\displaystyle Q^+ _{RR} (g,f) &=& \displaystyle\int _{\R^N} \, dv_*
 \int _{\ens{S}^{N-1}} \, d\sigma \, \Phi_{R,n} \, b_{R,m} \, g'_* \, f'.
\end{array}\end{cases}
 \end{equation*}

%%%%%%%%%%%%%%%%%%%%%%%%%%%%%%%%%%%%%%%%%%%%%%%%%%%%
%%%%%%%%%%%%%%%%%%%%%%%%%%%%%%%%%%%%%%%%%%%%%%%%%%%%

\subsection{Regularity and integrability for smooth collision kernel}

In this section we shall prove the regularity property of the gain
operator under the assumption that both $\Phi$ and $b$ are
smooth and compactly supported:
 \begin{equation}\label{eq:hyplions}
 \Phi \in C^\infty _0 (\R^N \setminus \{0\}), \qquad
 b \in C^\infty _0 (-1,1).
 \end{equation}
The assumption~\eqref{eq:hyplions} is obviously satisfied by
the smooth part $Q^+ _S$ of the gain operator in the decomposition above.
Thus, the results in this section will apply to the mollified
operator $Q^+ _S$ in~\eqref{eq:splitQ^+S}.
Our main result in this section is the
 \begin{theorem}\label{theo:regQ^+S}
 Let $B(v-v_*,\sigma) = \Phi(|v-v_*|) b(\cos \theta)$ satisfy the
 assumption~\eqref{eq:hyplions}.
 Then, for all $s \in \R^+$, $\eta \in \R$,
  \begin{eqnarray}\label{eq:regQ^+S}
  \left\|Q^+(g,f)\right\|_{H^{s+\frac{N-1}{2}} _\eta} &\le&
  \Creg (s,B) \left\|g\right\|_{H^s _\eta} \left\|f\right\|_{L^1 _{2|\eta|}} \\ \nonumber
  \left\|Q^+(f,g)\right\|_{H^{s+\frac{N-1}{2}} _\eta} &\le&
  \Creg (s,B) \left\|g\right\|_{H^s _\eta} \left\|f\right\|_{L^1 _{2|\eta|}}
  \end{eqnarray}
 where the constant $\Creg (s,B)$ only
 depends on $s$ and on the collision kernel, see
 formulas~\eqref{eq:Creg} and~\eqref{eq:comput:Creg}.
 \end{theorem}

\Remark Of course assumption~\eqref{eq:hyplions} is left invariant
under the change $\theta\to \pi/2-\theta$, and therefore the
estimates in~\eqref{eq:regQ^+S} are symmetric under exchange of
$f$ and $g$.

\begin{proof}[Proof of Theorem \ref{theo:regQ^+S}]
We shall proceed in four steps, following the method
of Wennberg~\cite{Wenn:rado:94}. We shall make use of the elementary
lemma~\ref{lem:OPD} in the Appendix to explicitly control an error
term disregarded in~\cite{Wenn:rado:94}.
\medskip

%%%%%%%%%%% Etape 1 de la preuve de regularisation de Q^+
\noindent{\bf Step 1: The Carleman representation}
\smallskip

The idea of Carleman representation (see \cite{Carl:fond:32,Carl:fond:57}) is to parametrize
$Q^+$ by the variables $v'$ and $v' _*$ instead of $v_*$ and $\sigma$. This
change of variable leads to
 \begin{equation*}\label{eq:carl}
 Q^+(g,f) = \int_{\R^N} dv' \int_{E_{v,v'}} dv' _* \, \frac{\Phi(|v-v_*|) \,
 b(\cos \theta)}{|v-v'|^{N-1}} g' _* f'
 \end{equation*}
where $E_{v,v'}$ denotes for the hyperplan orthogonal to $v-v'$ and
containing $v$. Since $|v-v'|/|v-v_*| = \sin (\theta/2)$, we can parametrize
the kernel by
 \begin{equation*}\label{kernelcarl}
 \frac{\Phi(|v-v_*|) \, b(\cos \theta)}{|v-v'|^{N-1}}
 = \cal{B}(v-v_*,|v-v'|)
 \end{equation*}
where
 \begin{equation*}
 \cal{B}(v_1,|v_2|) = \frac{\Phi(|v_1|) \, 
 b\left(1-2\left(\frac{|v_2|}{|v_1|}\right)^2 \right)}{|v_2| ^{N-1}}.
 \end{equation*}
The fact that $\cal{B}$ is radial according the first variable
will not be used in the next step, but will reveal useful in Step~3 where 
some modified versions of the collision kernel will be needed.

Following~\cite{Wenn:rado:94}, we define,
for $w \in \ens{S}^{N-1}$ and $r,s \in \R$,
 \begin{equation*}\label{eq:transfradon}
 R_{w,r}g(s) =  \int _{w^\bot} \cal{B} (z +sw,r)\, g(z+sw)\,dz
 \end{equation*}
where $w^\bot$ denotes the hyperplane orthogonal to $w$ going
through the origin (this is a weighted Radon transform).
Then, for $y\neq 0$ we set
 \begin{eqnarray*}\label{eq:T}
 Tg(y) &=& \left[ R_{y/|y|,|y|} \right] g(|y|) \\ \nonumber
 &=& \int _{y + y^\bot} \cal{B} (z,y)\, g(z)\,dz.
 \end{eqnarray*}
By an easy computation,
 \begin{equation*}\label{eq:lien}
 Q^+ (g,f) = \int _{\R^N} f(v')\,(\tau _{v'}
 \circ T \circ \tau _{-v'}) \, g(v) \,dv'
 \end{equation*}
(this is the last formula in~\cite[section~2]{Wenn:rado:94}). Thus
it becomes clear that regularity estimates on the Radon transform
$T$ will result in regularity estimates on $Q^+$. More precisely,
a careful use of Fubini and Jensen theorems leads to
 \begin{equation*}
 \left\|Q^+ (g,f)\right\|^2 _{H^{s+\frac{N-1}{2}} _\eta (\R^N)}
 \le \left\|f\right\|_{L^1} \int _{\R^N} |f(v')|
 \Bigl \|(\tau _{-v'} \circ T \circ \tau _{v'})g(v)
 \Bigr \|^2 _{H^{s+\frac{N-1}{2}} _\eta (\R^N)} dv',
 \end{equation*}
and we see that
 \begin{equation*}
 \left\|Q^+ (g,f)\right\| _{H^{s+\frac{N-1}{2}} _\eta (\R^N)}
 \le \Creg (s,B) \left\|f\right\|_{L^1 _{2|\eta|}}
\left\|g\right\|_{H^s _\eta (\R^N)},
 \end{equation*}
if we define $\Creg (s,B)$
as the best constant in the inequality
 \begin{equation} \label{eq:Tgineq}
 \left\|Tg\right\|_{H^{s+\frac{N-1}{2}} _\eta (\R^N)}
 \le C \left\|g\right\|_{H^s _\eta (\R^N)}.
 \end{equation}
\medskip

%%%%%%%%%% Etape 2 de la preuve de regularisation de Q^+
\noindent{\bf Step 2: Estimates of radial derivatives of $T$}
\smallskip

We now start to establish~\eqref{eq:Tgineq}. As we shall see in
the next step, it suffices to study the regularity with respect to
the modulus of the relative velocity variable, because the angular
derivatives can be controlled by the radial ones. We shall work in
spherical coordinates and write $Tg(r w) = R_{w,r} g(r)$ ($r > 0$,
$w \in \ens{S}^{N-1}$). We introduce the ``radial Fourier
transform'', $\RF$, and the Fourier transform in $\R^N$, $\F$, by
the formulas
\[ \RF f (\rho w) = \frac{1}{(2\pi)^{1/2}} \int_{\R} dr\,
 e^{i\rho r} f(r w), \]
\[ \F f(\xi) = \frac{1}{(2\pi)^{N/2}} \int_{\R^N} dv\,
e^{iv\cdot \xi} f(v). \]

In particular,
 \begin{equation*}
 \RF \left[ \langle r \rangle^\eta Tg \right](\rho w) =
 \frac{1}{(2\pi)^{1/2}} \int_{\R} dr\,
 e^{i\rho r} \langle r \rangle^\eta \int_{w^\bot} dz \,
 \cal{B} (z +rw,r) g(z+rw).
 \end{equation*}
Let $u=z + r w$. By Fubini's theorem and some simple computations,
 \begin{equation*}
 \RF \left[ \langle r \rangle^\eta Tg \right](\rho w)
 = (2\pi)^{\frac{N-1}{2}}
 \F \big[ g(\cdot)
 \cal{B}(\cdot,|(\cdot,w)|) \langle (\cdot,w) \rangle^\eta \big] (\rho w).
 \end{equation*}
By this we can estimate the $H^{s+\frac{N-1}{2}} _\eta$ norm according to the
radial variable. Let
 \[ \left\|Tg\right\|^2 _{H^{s+\frac{N-1}{2}} _\eta (\ens{S}^{N-1} \times \R)}
 = \int_{\ens{S}^{N-1}} dw \int_{\R} d\rho \,
 \langle \rho \rangle^{2\left(s+\frac{N-1}{2}\right)}
 \Bigl |\RF \big[ \langle r \rangle^\eta Tg \big] (\rho w)
 \Bigr |^2 \]
 \[ = (2\pi)^{N-1} \int_{\ens{S}^{N-1}} dw
 \int_{\R} d\rho \, \langle \rho \rangle^{2\left(s+\frac{N-1}{2}\right)}
 \Bigl |\F \big[ g(\cdot)
 \cal{B}(\cdot,|(\cdot,w)|) \langle (\cdot,w)
 \rangle^\eta \big] (\rho w)\Bigr |^2. \]
We change variables to get back to Euclidean coordinates, and find
 \begin{eqnarray*}
 && \left\|Tg\right\|^2 _{H^{s+\frac{N-1}{2}} _\eta(\ens{S}^{N-1} \times \R)} \\ \nonumber
 &=& (2\pi)^{N-1} \int_{\R^N}
 \langle \xi \rangle^{2s+N-1} |\xi|^{-(N-1)}\left| \F
 \left[ g(\cdot)
 \cal{B}\left(\cdot,\left|\left(\cdot,\frac{\xi}{|\xi|}\right)\right|\right)
 \left\langle (\cdot,\frac{\xi}{|\xi|}) \right\rangle^\eta
 \right] (\xi)  \right|^2 d\xi.
 \end{eqnarray*}
Now we cut this expression into two parts: for $|\xi| > 1$, the
inequality $|\xi|^2 > 1/2 (1+ |\xi|^2)$ implies that the
right-hand side is bounded from above by
 \begin{multline*}
 (8\pi)^{N-1} \int _{|\xi| > 1} \langle \xi \rangle^{2s}
 \left|\F \left[ g(\cdot)
 \cal{B}\left(\cdot,\left|\left(\cdot,\frac{\xi}{|\xi|}\right)\right|\right)
 \left \langle (\cdot,\frac{\xi}{|\xi|}) \right\rangle^\eta \right]
 (\xi)\right|^2 \\
 + (2\pi)^{N-1} 2^{s+\frac{N-1}{2}}
 \left( \int_{\ens{B}^N} \frac{d\xi}{|\xi|^{N-1}}\right) \sup _{|\xi| \le 1}
 \left|\F \left[g(\cdot)
 \cal{B}\left(\cdot,
 \left| \left(\cdot,\frac{\xi}{|\xi|} \right)\right| \right )
 \left\langle (\cdot,\frac{\xi}{|\xi|})
 \right\rangle^\eta \right] (\xi)\right|^2,
 \end{multline*}
where $\ens{B}^N$ stands for the ball of radius~1.

Then, on one hand Lemma~\ref{lem:OPD} implies
 \begin{multline*}
 \int _{|\xi| > 1} \langle \xi \rangle^{2s} \left |\F \left[ g(\cdot)
 \cal{B}\left(\cdot,\left|\left(\cdot,\frac{\xi}{|\xi|}\right)\right|
 \right)
 \left\langle (\cdot,\frac{\xi}{|\xi|}) \right\rangle^\eta \right]
 (\xi) \right |^2 \\
 \le \left\|g\right\|_{H^s _\eta}^2
 \left\|\cal{B} \left (x, \left|\left(x, \frac{y}{|y|} \right)\right|
 \right) \frac{\left\langle (x,\frac{y}{|y|})
 \right\rangle^\eta}
 {\langle x \rangle^\eta}\right\|_{L^\infty _y (H^S _x)} ^2
 \end{multline*}
where $S = s + \lfloor N/2 \rfloor +1$. On the other hand,
for each $|\xi|\le 1$,
 \[ \left |\F \left[
 g(\cdot)\cal{B}(\ldots)
 \left\langle (\cdot,\frac{\xi}{|\xi|}) \right\rangle^\eta \right]
 (\xi)\right |
 = \left| \frac{1}{(2\pi)^{N/2}} \int _{\R^N} e^{-i \xi \cdot x} g(x)
 \cal{B}(\ldots) \left\langle (x,\frac{\xi}{|\xi|})
 \right\rangle^\eta dx \right|. \]
Hence, by the Cauchy-Schwarz inequality,
 \[ \left |\F \left[
 g(\cdot)\cal{B}(\ldots)
 \left\langle (\cdot,\frac{\xi}{|\xi|})
 \right\rangle^\eta \right] (\xi) \right |
 \le \frac{1}{(2\pi)^{N/2}} \left\|g\right\|_{L^2 _\eta}
 \sup_{w \in \ens{S}^{N-1}}
 \left\|\cal{B}(x,|(x,w)|) \frac{\left\langle (x,w)
 \right\rangle^\eta} {\langle x \rangle^\eta}
 \right\|_{L^2 (\R^N _x)}.\]
Adding up the previous inequalities, we conclude that
 \begin{equation}\label{eq:estim:step2}
 \left\|Tg\right\| _{H^{s+\frac{N-1}{2}} _\eta
 (\ens{S}^{N-1} \times \R)} \le
 \cst(N,s) \, \left\|g\right\|_{H^s _\eta}
 \sup_{w \in \ens{S}^{N-1}}
 \left\|\cal{B}(x,|(x,w)|)\frac{\langle (x,w) \rangle^\eta}
 {\langle x \rangle^\eta}\right\|_{H^S (\R^N _x)}.
 \end{equation}
\medskip

%%%%%%%%%%%%% Etape 3 de la preuve de regularisation de Q^+
\noindent{\bf Step 3: Corollary: estimates of the angular derivatives of $T$}
\smallskip

Here we show how to get estimates on the \emph{angular} derivatives of
$Tg$ thanks to the estimates on the radial derivatives.
We first require the exponent $s+\frac{N-1}{2}$ to be
integer, so that the $H^{s+\frac{N-1}{2}} _\eta$ norm can be
computed in terms of norms of derivatives. Then,
 \begin{equation*}
 \frac{\partial (Tg)}{\partial y_i} (y) = \sum_j
 \frac{\partial w_j(y)}{\partial y_i} \frac{\partial}{\partial w_j} R_{w,s}g(s)
 + \frac{\partial s(y)}{\partial y_i} \left[\frac{\partial}{\partial s} R_{w,r}g(s)\right]_{s=r}
 + \frac{\partial r(y)}{\partial y_i} \left[\frac{\partial}{\partial r} R_{w,r}g(s)\right]_{s=r}
 \end{equation*}
where
 \[ w(y) = \frac{y}{|y|}, \qquad r(y) = s(y) = |y|, \]
and higher-order variants of this formula can obviously be obtained by
differentiating at arbitrary order. Let us assume that $\supp(\Phi) \subset
[\alpha,+\infty)$ and $\supp(b) \subset [\var,1]$ ($\alpha >0$ and $0<\var<1$).
Then $\supp(\mathcal{B}) \subset [\alpha,+\infty) \times [\var \alpha ,+\infty)$,
and one can easily establish that
 \[ \left| \frac{\partial^\nu w_j(y)}{\partial y^\nu} \right| \leq
 \frac{\cst(N)}{(\alpha \var)^{|\nu|}}, \qquad \left|
 \frac{\partial^\nu s(y)}{\partial y^\nu} \right|, \quad \left|
 \frac{\partial^\nu r(y)}{\partial y^\nu} \right| \leq
 \frac{\cst(N)}{(\alpha \var)^{|\nu|-1}}\]
in the support of ${\mathcal B}$.

Our second tool is the following property of the Radon transform: it
can be rewritten
 \begin{eqnarray*}
 R_{w,r} g (s) &=& \int _{w^{\bot}} g(z+sw) \mathcal{B} (z+sw,r) \, dz \\
               &=& \int _{\R^N} g(u) \mathcal{B} (u,r) \delta (w \cdot u - s ) \, du.
 \end{eqnarray*}
where $\delta$ is the Dirac mass at $0$ on $\R$. Thus
 \begin{eqnarray*}
 \frac{\partial}{\partial w_j} R_{w,r} g (s) &=&
 \int _{\R^N} g(u) \mathcal{B} (u,r) u_j\,
 \delta' (w \cdot u - s ) \, du \\ 
 &=& -\,\frac{\partial}{\partial s} \int _{\R^N} g(u)
 \mathcal{B} (u,r) u_j\, \delta (w \cdot u - s ) \, du \\ 
 &=& -\,\frac{\partial}{\partial s} \widetilde{R}_{w,r} (g) (s)
 \end{eqnarray*}
where $\widetilde{R}_{w,r}$ is defined by the new kernel
$\mathcal{B} (u,r) u_j$.
Thus the angular derivative
$\frac{\partial}{\partial w_j} R_{w,r} g (s)$ can be obtained from
the estimate~\eqref{eq:estim:step2} of Step~2, upon changing the collision kernel
by another one, only differing by a factor of $u_j$. The same
holds true for all order derivatives...

To conclude with the regularization property of $T$, it is enough
to notice that the derivatives along $r$ are already taken into
account above, and to use the above-mentioned
commutation property for the angular derivatives. We conclude that
equation~\eqref{eq:Tgineq} holds true with
 \begin{eqnarray}\label{eq:Creg}
 \Creg (s,B) &=& \frac{\cst(s,N)}{(\alpha \var)^{s+\frac{N-1}{2}}} \\ \nonumber
 && \sup \left\{
 \left\| \cal{B}(x,|(x,w)|) x^\nu
 \frac{\langle (x,w) \rangle ^\eta}{\langle x \rangle ^\eta} \right\|
 _{H^{S-|\eta|} (\R^N _x)}; \
 |\nu|\le s+\frac{N-1}{2}; \quad w \in \ens{S}^{N-1} \right\}.
 \end{eqnarray}
 \medskip

This concludes the proof of~\eqref{eq:Tgineq} when $s+(N-1)/2$ is
an integer. The general case follows by the Riesz-Thorin interpolation
theorem again.
\end{proof}

\noindent{\bf Order of the constant according to the convolution
parameters:}
\smallskip

The computation of an upper bound on the constant $\Creg (s,B)$
for the collision kernel $\Phi_{S,n} b_{S,m}$ according to the
mollifying parameters $m$ and $n$ is tedious but straightforward.
One may easily obtain a polynomial bound in the form
 \begin{eqnarray} \nonumber
 \Creg(s,B) &\le& \cst(s,N) m^{as+b} n^{a's+b'}
 \left\|1_{\ens{I}_m} b\right\|_{L^1 (\ens{S}^{N-1})} \\ \label{eq:comput:Creg}
 &\leq& \cst(s,N) m^{as+b} n^{a's+b'} \left\| b\right\|_{L^1 (\ens{S}^{N-1})}.
 \end{eqnarray}
where $a,a',b,b'$ stand for some constant depending only the
dimension $N$ and $\gamma$.
\medskip

We conclude this section with the following corollary of
Theorem~\ref{theo:regQ^+S}, which translates the gain of
regularity into a gain of integrability.

 \begin{corollary}\label{coro:gain:leb}
 Let us consider a collision kernel $B$ satisfying the smoothness
 assumption~\eqref{eq:hyplions}.
 Then, for all $p \in (1;+\infty)$, $\eta \in \R$, we have
  \begin{eqnarray*}\label{eq:intQ^+S}
  \left\|Q^+(g,f)\right\|_{L^q _\eta} &\le&
  \Cint (p,\eta,B) \left\|g\right\|_{L^p _\eta} \left\|f\right\|_{L^1 _{2|\eta|}} \\ \nonumber
  \left\|Q^+(f,g)\right\|_{L^q _\eta} &\le&
  \Cint (p,\eta, B) \left\|g\right\|_{L^p _\eta} \left\|f\right\|_{L^1 _{2|\eta|}}
  \end{eqnarray*}
 where the constant $\Cint (p,\eta,B)$ only depends on
 the collision kernel, $p$ and $\eta$, and $q>p$ is given by
  \begin{equation*}
  q =
  \left\{
  \begin{array}{cl} \displaystyle
  \frac{p}{2-\frac{1}{N} + p \left(\frac{1}{N} -1 \right)}
  & \mbox{if } p \in (1;2] \\ \\
  pN & \mbox{if } p \in [2;+\infty).
  \end{array}
  \right.
  \end{equation*}
 \end{corollary}

\Remark Just as $\Creg$, the constant $\Cint (m,n)$ depends on 
the mollifying parameters in a polynomial way. Note that 
the constant $\Cint (p,\eta, B)$ in Corollary~\ref{coro:gain:leb}
does not depend on the weight exponent $\eta$ anymore in the quadratic case (just as
in section~\ref{sec:prelim}).

\begin{proof}
The proof is almost obvious. When $p=2$, it is a direct consequence of
Theorem~\ref{theo:regQ^+S} with $s=0$, and the Sobolev injection
$H^{\frac{N-1}{2}} _\eta \hookrightarrow L^{2N} _\eta$ (with a
constant only depending on $N$). The general case follows by a
Riesz-Thorin interpolation betweeen this estimate and the
convolution-like inequalities in Theorem~\ref{theo:conv}.
\end{proof}

%%%%%%%%%%%%%%%%%%%%%%%%%%%%%%%%%%%%%%%%%%%%%%%%%%%%
%%%%%%%%%%%%%%%%%%%%%%%%%%%%%%%%%%%%%%%%%%%%%%%%%%%%

\subsection{Regularity and integrability for nonsmooth collision kernel}

In this paragraph we extend the regularity of $Q^+$ to general
nonsmooth kernels. There are at least two strategies for that,
which will lead to slightly different results. We shall first
give a general result of ``gain of integrability/regularity'', in
a form which is remindful of the classical Povzner inequalities
used to study the $L^1$ moment behavior (besides it will play the same
role in the proof of propagation of $L^p$ moments).
\medskip

\bul The following inequality will turn
out to be the most appropriate for our study of propagation of
integrability. We state it only in its quadratic version, the
bilinear version would be slightly more intricate but easy to
write down as well.

 \begin{theorem} \label{theo:gain:int}
 Let $B$ be a collision kernel satisfying
 assumptions~\eqref{eq:hypcut-off},~\eqref{eq:hyprad}
 and~\eqref{eq:hypang2}. Then, for all $p > 1$,
 $k>\gamma$ and $\eta \geq
 -\gamma$, there exist constants $C$ and $\kappa$, and
 $q<p$ ($q$ only depending on $p$ and $N$), such that for all
 $\var>0$, and for all measurable $f$,
  \begin{equation*} \label{inegintQ}
  \|Q^+(f,f)\|_{L^p_\eta} \leq C \var^{-\kappa}\,
  \|f\|_{L^q_\eta}\|f\|_{L^1_{2|\eta|}} + \var \|f\|_{L^p_{\gamma
  +\eta}} \|f\|_{L^1_{k+2\eta_+}}.
  \end{equation*}
 \end{theorem}
This estimate expresses a ``mixing'' property of the
$Q^+$ operator: the dominant norm $L^p_{\gamma+\eta}$ appears with
a constant $\var$ as small as desired; and for the rest, we can
lower both the Lebesgue exponent and its weight. This property is of 
course consistent with the compactness properties of $Q^+$, and
in complete contrast with the properties of the loss term $Q^-$.

\begin{proof}[Proof of Theorem~\ref{theo:gain:int}]
We split $Q^+$ as $Q^+ _{S} + Q^+ _{RS} + Q^+_{SR} +
Q^+ _{RR}$ and we shall estimate each term separately. 
From the beginning we assume, without loss of generality, that the
angular kernel $b(\cos\theta)$ has support in $[0,\pi/2]$.
Remember that the truncation parameters $n$ (for the kinetic part)
and $m$ (for the angular part) are implicit in the decomposition
of $Q^+$.

By Corollary~\ref{coro:gain:leb}, there exists a constant
$\Cint(m,n)$, blowing up polynomially as $m\to\infty$,
$n\to\infty$, such that
 \[ \left\|Q^+ _{S} (f,f)\right\| _{L^p_\eta}
 \leq \Cint(m,n) \|f\|_{L^q_\eta} \|f\|_{L^1_{2|\eta|}}, \] for
some $q<p$, namely
 \begin{equation*}\label{eq:defq}
 q= \begin{cases}
 \displaystyle \frac{(2N-1)p}{N+(N-1)p} \quad \text{if $p \in (1;2N]$} \\ \\
 \displaystyle \frac{p}{N} \quad \text{if $p \in [2N;+\infty)$}
 \end{cases}
 \end{equation*}
(the roles of $p$ and $q$ are exchanged here with respect to
Corollary~\ref{coro:gain:leb}....)

Next, we shall take advantage of the fact that $b_{R,m}$ has a
very small mass (assumptions~\eqref{eq:hypcut-off}
and~\eqref{eq:hypang2}), and write, using
Corollary~\ref{coro:conv} with $k=\gamma$,
 \[ \|Q^+_{RR} (f,f) \|_{L^p_\eta} \leq C m^{-\delta}
 \|f\|_{L^1_{|\gamma+\eta|+|\eta|}} \|f\|_{L^p_{\gamma+\eta}}, \]
for some constant $C$ only
depending on $C_{\Phi}$.
A similar estimate holds true for $\|Q^+_{SR}\|_{L^p_\eta}$. Since
$\gamma+\eta\geq 0$, we can write $|\gamma+\eta|+|\eta| =
\gamma+2\eta_+$, where $\eta_+= \max(\eta,0)$.

It remains to estimate the term $Q^+_{RS}$. For this we shall
consider separately large and small velocities, and write
$f=f_r+f_{r^c}$, where
 \begin{equation*}
 \left\{
 \begin{array}{l}
 f_r = f\, 1_{\{|v|\le r\}} \\
 f_{r^c} = f\, 1_{\{|v| > r\}}.
 \end{array}
 \right.
 \end{equation*}
On the one hand, we use Theorem~\ref{theo:conv}, and pick a
$k>\gamma$, in order to ensure that
$\left\|\Phi_{R,n}\right\|_{L^\infty _{-k}}$ goes to~0 as
$n\to\infty$. Thanks to the H\"older assumption~\eqref{eq:hyprad},
one can easily prove
\[ \left\|\Phi_{R,n}\right\|_{L^\infty _{-k}} \le \cst \
 n^{-\left(\min(\gamma,k-\gamma)\right)}.\]

It follows
 \begin{align*}
 \left\|Q^{+} _{RS}(f,f_r) \right\| _{L^p _\eta}
 & \leq C \left\|f \right\| _{L^1 _{|k+\eta|+|\eta|}}
 \left\|f_r \right\|^{p} _{L^p _{k+\eta}} \|\Phi_{R,n}\|_{L^\infty_{-k}}\\
 & \leq C \left\|f \right\| _{L^1 _{|k+\eta|+|\eta|}} r^{k-\gamma}
 \|f\|_{L^p_{\gamma+\eta}} n^{\gamma-k}\\
 & \leq C \left ( \frac{r}{n} \right )^{k-\gamma}
 \|f\|_{L^1_{k+2\eta_+}} \|f\|_{L^p_{\gamma+\eta}}.
 \end{align*}
(here $\theta_b = \pi/2$ thanks to the symmetrization).

\Remark This is the only place where we use a regularity estimate
on $\Phi$.
\smallskip

On the other hand, the support of $b_{S,m}$ lies a positive
distance ($O(1/m)$) away from~0, so~\eqref{eq:hypangconvol2} holds
true with $\theta_b =\cst\, m^{-1}$. Thus we can apply
Theorem~\ref{theo:conv} with $f$ and $g$ exchanged, to find
 \[
 \left\|Q^{+} _{RS}(f,f_{r^c}) \right\| _{L^p _\eta}
 \leq C m^\beta \|f_{r^c}\|_{L^1_{|\gamma+\eta|+|\eta|}}
 \|f\|_{L^p_{\gamma+\eta}}.\]
where $\beta = \max(-\eta,0)+2/p'$ and $C$ depends only on $C_{\Phi}$.
Since we assume $\gamma+\eta \geq 0$, this can also be bounded by
 \[ C m^\beta r^{\gamma-k} \|f\|_{L^1_{k+\eta+|\eta|}}
 \|f\|_{L^p_{\gamma+\eta}} = C m^\beta r^{\gamma-k}
 \|f\|_{L^1_{k+2\eta_+}} \|f\|_{L^p_{\gamma+\eta}}. \]

To sum up, we have obtained
 \[ \|Q^+(f,f) \|_{L^p_\eta} \leq
 C_1(m,n) \|f\|_{L^q_\eta} \|f\|_{L^1_{2|\eta|}} + C \Bigl [
 m^{-\delta} +  \left ( \frac{r}{n} \right )^{k-\gamma} +
 \frac{m^\beta}{r^{k-\gamma}} \Bigr ] \|f\|_{L^1_{k+2\eta_+}}
 \|f\|_{L^p_{\gamma+\eta}}. \]
The conclusion follows by choosing first $m$ large enough, then
$r$, then $n$.
\end{proof}

We turn to another similar theorem in which the emphasis is laid
on regularity rather than integrability and whose proof is quite
similar.

 \begin{theorem} \label{theo:gain:reg}
 Let $B$ be a collision kernel satisfying
 assumptions~\eqref{eq:hypcut-off},~\eqref{eq:hyprad}
 and~\eqref{eq:hypang2}. Then, for all $s > 0$, $k>\gamma$ and $\eta\geq
 -\gamma$, there exist constants $C$ and $\kappa$, and $0 \leq s' < s$
 ($s' = \max\left(s-\frac{N-1}{2},0\right)$ only depending on $s$ and $N$),
 such that for all $\var$, and for all measurable $f$,
  \begin{equation*} \label{inegintQbis}
  \|Q^+(f,f)\|_{H^s_\eta} \leq C \var^{-\kappa}\, \|f\|_{H^{s'}
  _\eta}\|f\|_{L^1_{2|\eta|}} + \var \|f\|_{H^s_{\gamma +\eta}}
  \|f\|_{W^{\lceil s \rceil,1} _{k+2\eta_+}}.
  \end{equation*}
 \end{theorem}

\begin{proof}[Proof of theorem~\ref{theo:gain:reg}]
The proof follows the same path as the previous one. The term $Q^+
_S$ is estimated by Theorem~\ref{theo:regQ^+S}, the terms $Q^+ _{SR}$ and
$Q^+ _{RR}$ are estimated by Theorem~\ref{theo:conv}. For the
remaining term $Q^+ _{RS}$, we also estimate separately large and
small velocities. But this time, the splitting $f=f_r+f_{r^c}$
should be
 \begin{equation*}
 \left\{
 \begin{array}{l}
 f_r = f\, \chi _r \\
 f_{r^c} = f - f_r.
 \end{array}
 \right.
 \end{equation*}
where $\chi _r$ is a $C^\infty$ function with bounded derivatives
and such that $\chi _r = 1$ on $|v| \le r$ and $\supp \chi _r
\subset B(0,r+1)$. The end of the proof is straightforward.
\end{proof}
%\Remark By using Theorem~\ref{theo:reg:Q^+:BD} below and minor
%changes in the proof above, one could also prove,
%under the additional assumption~\eqref{eq:hypL^2} below, a functional
%inequality of the type
% \begin{equation} \label{inegintQ2}
% \|Q^+(f,f)\|_{H^s_\eta} \leq C \var^{-\kappa}\, \|f\|_{H^{s'}
% _\eta}\|f\|_{L^1_{2|\eta|}} + \var \|f\|_{H^{s'}_{\gamma +\eta}}
% \|f\|_{H^{s'} _{\eta+\gamma+1}}.
% \end{equation}
% CM : j'enleve cette remarque que j'avais mise et qui me parait un peu douteuse...
Note that there are other possible variants as well....
\smallskip
% {\bf Je ne suis pas sur que ce soit le theoreme qu'il nous faille
% ! Peut-etre faut-il se limiter a
%  \begin{equation} \label{inegintQ}
%  \|Q^+(f,f)\|_{L^2_\eta} \leq C \var^{-\kappa}\,
%  \|f\|_{H^{s} _\eta}\|f\|_{L^1_{2|\eta|}} + \var \|f\|_{L^2_{\gamma
%  +\eta}} \|f\|_{L^1_{k+2\eta_+}}
%  \end{equation}
% pour un certain $s>0$. }
% \end{theorem}

\bul The first way to a regularity result for the full kernel is
to use the method by Bouchut and Desvillettes
in~\cite{BoucDesv:rgain:98}. Hence it is possible to extend
Theorem~2.1 in~\cite{BoucDesv:rgain:98} into the following
 \begin{theorem}\label{theo:reg:Q^+:BD}
 Let $B(v-v_*,\sigma) = \Phi(|v-v_*|) \, b (\cos \sigma)$ be a
 collision kernel such that $\Phi$ satisfies the
 assumption~\eqref{eq:hyprad} and $b$ satisfies
  \begin{equation} \label{eq:hypL^2}
  \left\|b\right\|_{L^2
  (\ens{S}^{N-1})} < + \infty
  \end{equation}
 in the sense that $\int b(\cos\theta)^2\,\sin^{N-1}\theta\,d\theta
 <+\infty$. Then for all $s \geq 0$ and $\eta \geq 0$,
  \begin{equation*}\label{eq:reg:Q^+:BD}
  \left\| Q^+ (g,f) \right\|_{H^{s + \frac{N-1}{2}} _{\eta}}
  \le \Cbd \left[ \left\|g\right\|_{H^s _{\eta+ \gamma +1}} \left\|f\right\|_{H^s
  _{\eta+\gamma+1}} + \left\|g\right\|_{L^1 _{\eta + \gamma}} \left\|f\right\|_{L^1 _{\eta + \gamma}} \right]
  \end{equation*}
 where $\Cbd$ only depends on $N$ and on $\left\|b\right\|_{L^2
 (\ens{S}^{N-1})}$.
 \end{theorem}

\Remarks 1. Of course assumption~\eqref{eq:hypL^2} is stronger
than~\eqref{eq:hypcut-off}; it is however still reasonable in the
context of cut-off hard potentials (in particular for hard
spheres, in which $b$ is just a constant).
\smallskip

2. The inequality here is not adapted to our study of
integrability, but will be useful for our study of regularity.
Moreover, the proof and the constants are simpler than those which
led us to Theorem~\ref{theo:reg:Q^+:full}.
\medskip

\bul The second way towards a regularity result for the full
kernel is to combine Theorem~\ref{theo:conv} and
Theorem~\ref{theo:regQ^+S} and make an explicit interpolation. By
this one can prove the

% \begin{theorem}\label{theo:reg:Q^+:full}
% Let us consider a collision kernel $B$ satisfying
% assumptions~\eqref{eq:hypcut-off},~\eqref{eq:hyprad}, \eqref{eq:hypang2} and~\eqref{eq:hypangconvol}.
% Then for all $k>\gamma$
% there exists $\alpha > 0$, depending only on $B$,
% such that for all $s \ge 0$ and $\eta \in \R$
%  \begin{equation}
%  \left\|Q^+(g,f)\right\|_{H^{s+\alpha} _\eta}
%  \le C \left\|g\right\|_{W^{\lceil s \rceil,1} _{k+2\eta+}}
%  \left\|f\right\|_{H^s _{\gamma+\eta}},
%  \end{equation}
% for some constant $C$ which only depends on $s$ and $B$.
% \end{theorem}
%
% \begin{corollary}\label{theosym}
% Let us consider a collision kernel $B$ satisfying
% assumptions~\eqref{eq:hypcut-off},~\eqref{eq:hyprad} and~\eqref{eq:hypang2}.
% Then for all $k>\gamma$ there exists $\alpha > 0$, depending only on $B$,
% such that for all $s \ge 0$ and $\eta \in \R$
%  \begin{equation}\label{eq:reg:Q^+:full}
%  \left\|Q^+(f,f)\right\|_{H^{s+\alpha} _\eta}
%  \le C \left\|f\right\|_{W^{\lceil s \rceil,1} _{k+2\eta+}}
%  \left\|f\right\|_{H^s _{\gamma+\eta}},
%  \end{equation}
% for some constant $C$ which only depends on $s$ and $B$.
% \end{corollary}
 \begin{theorem}\label{theo:reg:Q^+:full}
 Let us consider a collision kernel $B$ satisfying
 assumptions~\eqref{eq:hypcut-off},~\eqref{eq:hyprad} and~\eqref{eq:hypang2}.
 Then for all $k>\gamma$ and $\eta\geq
 -\gamma$, there exists $\alpha > 0$, depending only on $B$,
 such that for all $s \ge 0$ and $\eta \in \R$
  \begin{equation*}\label{eq:reg:Q^+:full}
  \left\|Q^+(f,f)\right\|_{H^{s+\alpha} _\eta}
  \le C \left\|f\right\|_{W^{\lceil s \rceil,1} _{k+2\eta+}}
  \left\|f\right\|_{H^s _{\gamma+\eta}},
  \end{equation*}
 for some constant $C$ which only depends on $s$ and $B$.
 \end{theorem}
% CM : j'ai enleve la version bilineaire qui me parait fausse telle quelle 
%      (les estimations de petitesse sont tres liees au sens dans lequel on les fait)

\begin{proof}[Proof of~Theorem~\ref{theo:reg:Q^+:full}]
Let us take $s \in \R_+$ and $\eta \in \R$.
We have the following estimates on the
four parts of the decomposition of $Q^+$ (by symmetrization the angular part of the collision
kernel is supposed to be zero for $\theta \ge \pi/2$).
\medskip

{\bf 1. For the smooth part}, Theorem~\ref{theo:regQ^+S} gives
 \begin{equation*} \label{eq:estim:Q^+S}
 \left\|Q^+ _S (f,f)\right\| _{H^{s+\frac{N-1}{2}} _\eta (\R^N)}
 \le C_1 \left\|f\right\|_{L^1 _{2|\eta|}} \left\|f\right\|_{H^s _\eta (\R^N)}
 \end{equation*}
where $C_1 = \Creg (m,n)$ blows up polynomially as $m\to\infty$,
$n\to\infty$.
\medskip

{\bf 2. To control the effect of small deviation angles}, we use
again Corollary~\ref{coro:conv}, and the dependence of the constant on
$\left\|b_{R,n}\right\|$ to ensure it goes to zero; we obtain
as in the proof of Theorem~\ref{theo:gain:int}
 \begin{equation*} \label{eq:estim:Q^+R2}
 \left\|Q^+ _{SR} (f,f),\>  Q^+ _{RR} (f,f)\right\|_{H^s _\eta} \le
 C_2 \left\|f\right\|_{W^{\lceil s \rceil,1} _{\gamma+2\eta+}}
 \left\|f\right\|_{H^s _{\gamma+\eta}}
 \end{equation*}
where $C_2 = \cst(N) \ \left\|b_R \right\|_{L^1(\ens{S}^{N-1})}
\left\|\Phi\right\|_{L^\infty _{-\gamma}}$, which thanks to
assumption~\eqref{eq:hypang2} can be bounded
from above by $\cst(C_B,N) m^{-\delta}$.
\medskip

{\bf 3. To control the effect of singularities of the kinetic kernel and
high velocities},
we use again Theorem~\ref{theo:conv} and pick a $k>\gamma$. As in the
proof of Theorem~\ref{theo:gain:reg}, we prove
 \begin{equation*}\label{eq:estim:Q^+R1}
 \Bigl\|Q^+ _{RS} (f,f)\Bigr\|_{H^s _\eta} \le C_3
 \left\|f\right\|_{W^{\lceil s \rceil,1} _{k+2\eta+}}
 \left\|f\right\|_{H^s _{\gamma+\eta}}
 \end{equation*}
where $C_3 = C \Bigl [ m^{-\delta} +  \left ( \frac{r}{n}
\right)^{k-\gamma} + \frac{m^\beta}{r^{k-\gamma}} \Bigr ]$, which
goes to $0$ polynomially according to the parameter $m$ when one
set $r$ then $n$ as well-chosen functions of $m$.
\medskip

To sum up, we know that for all $m \ge 1$, one can decompose $Q^+$
as $Q^+ = Q^+ _{S,m} + Q^+ _{R,m}$ (remember $n$ is now set as a
function of $m$), with the estimates
 \begin{equation*} \label{eq:sum:interp:Q^+}
 \left\{
 \begin{array}{l}
 \left\|Q^+ _{S,m} (f,f)\right\| _{H^{s+\frac{N-1}{2}} _\eta}
 \le C_1 \left\|f\right\|_{L^1 _{2|\eta|}} \left\|f\right\|_{H^s _\eta} \\ \\
 \left\|Q^+ _{R,m} (f,f)\right\|_{H^s _\eta} \le
 (C_2 + C_3)
 \left\|f\right\|_{W^{\lceil s \rceil,1} _{k+2\eta+}}
 \left\|f\right\|_{H^s _{\gamma+\eta}}
 \end{array}
 \right.
 \end{equation*}
By applying Theorem~\ref{theo:interp:sum:gen} in the
Appendix, we can conclude that
 \begin{equation*}
 \left\|Q^+ (f,f)\right\| _{H^{s+\alpha} _\eta}
 \le C \left\|f\right\|_{W^{\lceil s \rceil,1} _{k+2\eta+}}
 \left\|f\right\|_{H^s _\eta}
 \end{equation*}
for some $0< \alpha < \frac{N-1}{2}$ depending on the exponents of
polynomial control for each term. This concludes the proof.
\end{proof}

\Remark Some closely related results can be
found in Wennberg~\cite{Wenn:rado:94}, the goal is however different:
in this reference the author searches for sufficient conditions
on the collision kernel $B$, to ensure that the $H^{(N-1)/2}$ bound still
holds true. Here on the contrary we allow general collision kernels,
but, as a natural price to pay, the regularization which we obtain 
is in general strictly less than a gain of $(N-1)/2$ derivatives.
\medskip

%%%%%%%%%%%%%%%%%%%%%%%%%%%%%%%%%%%%%%%%%%%%%%%%%%%%
%%%%%%%%%%%%%%%%%%%%%%%%%%%%%%%%%%%%%%%%%%%%%%%%%%%%
%%%%%%%%%%%%%%%%%%%%%%%%%%%%%%%%%%%%%%%%%%%%%%%%%%%%

\section{Propagation of $L^p$ estimates} \label{sec:Lp}
\setcounter{equation}{0}
In this section we are interested in the propagation of $L^p$
integrability of the solutions of Boltzmann's equation and its
derivatives. Our proofs will be based on a differential
inequality approach. Most of the hard job has been done in the 
functional study of the previous section, and the proofs will
be much less technical now.

The bounds that we establish here will later serve as the first 
step for our study of propagation of regularity via a semigroup 
approach.

%%%%%%%%%%%%%%%%%%%%%%%%%%%%%%%%%%%%%%%%%%%%%%%%%%%%
%%%%%%%%%%%%%%%%%%%%%%%%%%%%%%%%%%%%%%%%%%%%%%%%%%%%
\subsection{Main result}

 \begin{theorem}\label{theo:prop:L^p}
 Let $B(v-v_*,\sigma) = \Phi(|v-v_*|) b(\cos \theta)$ satisfy
 assumptions~\eqref{eq:hypcut-off}, \eqref{eq:hyprad},
\eqref{eq:hypang2}, \eqref{eq:hyplbB},
 let $1<p<+\infty$ and let $f_0$ be a nonnegative function in
 $L^1_2 \cap L^p(\R^N)$. Then, the unique solution $f$ of the
 Boltzmann equation with initial datum $f_0$ satisfies the
 estimates
  \begin{equation} \label{eq:diffineqp}
  \frac{d\left\|f\right\|^p _{L^p}}{dt} \le
  C_+ \left\|f\right\|_{L^p} ^{p(1-\theta)} -
  K_- \left\|f\right\|_{L^p _{\gamma/p}} ^p
  \end{equation}
 for some constants $C_+,K_->0$, $\theta\in (0,1)$ which only depend
 on $p$, $N$, $B$, on upper bounds on $\|f\|_{L^1_2}$ and $H(f)$, and on a
 lower bound on $\|f\|_{L^1}$.

 In particular, there is an explicit constant $C_p(f_0)$, only
 depending on $B$, on an upper bound on $\|f_0\|_{L^1_2}+\|f_0\|_{L^p}$,
 and on a lower bound on $\|f_0\|_{L^1}$, such that
  \[ \forall t\geq 0, \qquad
  \|f(t,\cdot)\|_{L^p} \leq C_p(f_0). \]
 Moreover, for any $t >0$ and any $\eta >0$, we know that
 $f(t,\cdot) \in L^p_\eta(\R^N)$. More precisely, for any $t_0>0$,
  \[ \sup_{t\geq t_0} \|f(t,\cdot)\|_{L^p_\eta} < +\infty. \]
 Once again this bound can be computed in terms of
 $B$, an upper bound on $\|f_0\|_{L^1_2}+\|f_0\|_{L^p}$,
 a lower bound on $\|f_0\|_{L^1}$, and a lower bound
 on $t_0$.
 \end{theorem}

\begin{proof}[Proof of Theorem~\ref{theo:prop:L^p}]
Here we shall just be content with establishing the necessary a
priori estimates. The proof of the theorem follows from standard
approximation arguments, known results on the unique solvability
of the Boltzmann equation, with bounds in, say, weighted
$L^\infty$ if the initial datum also satisfies such bounds
(see the references indicated in the Introduction).

Let $f$ be a solution to the Boltzmann equation, supposed to be in
$C^1(\R_t,L^p)$. Also, since the solution is differentiable in $L^p$,
 \[ \frac{1}{p}\,
 \frac{d\left\|f\right\|^p _{L^p}}{dt} =
 \int f^{p-1} Q^+(f,f) \, dv - \int f^{p-1} Q^-(f,f)\,dv.\]
By Proposition~\ref{prop:minor},
 \begin{equation} \label{etape1}
 - \int f^{p-1} Q^- \, dv \leq - K \int f^p (1 + |v|)^\gamma \, dv
 \leq -K_0 \left\|f\right\|^p _{L^p _{\gamma/p}}.
 \end{equation}

On the other hand, by H\"older's inequality,
 \[ \int f^{p-1} Q^+ _S (f,f) \, dv
 \leq \left [ \int f^p \right]^{\frac{p-1}{p}} \left
 [\int (Q^+ _S)^p \right ]^{\frac1p} \]
 \[ = \|f\|_{L^p}^{p-1} \|Q^+ _S(f,f)\|_{L^p}. \]
and
 \[ \int f^{p-1} Q^+ _R (f,f) \, dv
 = \int \bigl (f \langle v\rangle^{\gamma/p} \bigr )^{p-1}
 \frac{Q^+}{\langle v\rangle^{\frac{\gamma}{p'}}} \leq \left [ \int
 (f \langle v\rangle^{\gamma/p})^p \right]^{\frac{p-1}{p}} \left
 [\int (Q^+ _R \langle v\rangle^{-\gamma/p'} )^p \right ]^{\frac1p} \]
 \[ = \|f\|_{L^p_{\gamma/p}}^{p-1} \|Q^+ _R(f,f)\|_{L^p_{-\gamma/p'}}. \]
By using the estimates on $Q^+ _S$ and $Q^+ _R$ proved in Theorem~\ref{theo:gain:int}
with $\eta=-\gamma/p'$, $k=2$ and $\var=K_0/(2\|f\|_{L^1_2})$,
we can find a constant $C$, depending on $\|f\|_{L^1 _2}$, such that
 \begin{equation*}
 \int f^{p-1} Q^+(f,f) \, dv \leq
 C \|f\|_{L^q}\|f\|_{L^1} \|f\|_{L^p}^{p-1}
 + \var\|f\|_{L^1_2} \|f\|_{L^p_{\gamma/p}} ^p
 \end{equation*}
where $q$ is defined by~\eqref{eq:defq}. Combining this with
elementary Lebesgue interpolation and the conservation of mass and
energy, we deduce that there exists a $\theta\in (0,1)$, only
depending on $N$ and $p$, and a constant $C_0$, only depending on
$N$, $p$, $B$ and $\|f_0\|_{L^1_2}$, such that
 \begin{align*}
 \int f^{p-1} Q^+(f,f) \, dv & \leq C_0
 \|f\|_{L^p}^{1-p\theta}\|f\|_{L^p}^{p-1}
 + \frac{K_0}{2} \|f\|_{L^p_{\gamma/p}}^p \\
 & \leq C_0 \|f\|_{L^p}^{p(1-\theta)}
 + \frac{K_0}{2} \|f\|_{L^p_{\gamma/p}}^p.
 \end{align*}
This together with~\eqref{etape1} concludes the proof of the
differential inequality~\eqref{eq:diffineqp} with $C_+ = C_0$ and 
$K_+ = K_0 /2$. 

From this differential inequality we see that the time-derivative
of $\|f(t,\cdot)\|_{L^p}^p$ is bounded by a constant, and
therefore $f(t,\cdot)$ lies in $L^p$ for all times. Moreover, if
$\|f(t,\cdot)\|_{L^p_{\gamma/p}}$ ever becomes greater than
$(C_+ / K_-)^{\frac{1}{p\theta}}$, it follows from~\eqref{eq:diffineqp} that
$(d/dt)\|f(t,\cdot)\|_{L^p} \leq 0$. Since $\|f\|_{L^p_{\gamma/p}}
\geq \|f\|_{L^p}$, we conclude that that
 \[ C_p(f_0) := \max \left[ \left\|f_0\right\|_{L^p} ; \
 \left(\frac{C_+}{K_-} \right) ^{\frac{1}{p\theta}} \right] \]
is a uniform upper bound for $\|f(t,\cdot)\|_{L^p}$.

Next, for all $\eta \geq 0$, a similar argument leads to the
a priori differential inequality
 \begin{equation}\label{eq:apriori:eqdif:weight}
 \frac{d\left\|f\right\|^p _{L^p _\eta}}{dt} \le
 C_+ \left\|f\right\|_{L^p _\eta} ^{p(1-\theta)} -
 K_- \left\|f\right\|_{L^p _{\eta + \gamma/p}} ^p.
 \end{equation}
where $C_+, K_-$ now depend on the entropy and on some
$\|f\|_{L^1_s}$ norm for $s$ large enough (depending on $\eta$).
We deduce that $\|f\|_{L^p_\eta}$ norms are propagated,
uniformly in time, if the initial datum possesses $L^1$ moments
of high enough order. Let $t_0>0$ be arbitrarily small;
for $t\geq t_0$, we know that
all the quantities $\|f(t,\cdot)\|_{L^1_s}$ are bounded,
uniformly in time, for all $s$, and these inequalities
therefore hold true with uniform constants as soon as $t\geq t_0$.

We next turn to the property of moment generation, i.e. the proof that
$L^p_\eta$ norms are automatically bounded for positive times.
These results are the analogue of the well-known results of
$L^1$ moment generation for hard potential with
cut-off (see for instance~\cite[Theorem~4.2]{Wenn:momt:97}).
Let $t_0>0$ be arbitrarily small. Integrating the
inequality~\eqref{eq:diffineqp} in time from~0 to $t_0$, we obtain
 \begin{equation*} \label{eq:app:mts}
 \int _{0} ^{t_0} \left\|f(s,\cdot)\right\|^p _{L^p _{\gamma/p}} \le
 \frac{C_+}{K_-} \int _0 ^{t_0} \left\|f(s,\cdot)\right\|^{p-\theta} _{L^p}
 + \frac{1}{K_-} \Bigl (
 \left\|f_0\right\|^p _{L^p} - \left\|f(t_0,\cdot)\right\|^p _{L^p}
 \Bigr ),
 \end{equation*}
which implies
 \begin{equation*}
 \int _0 ^{t_0} \left\|f(s,\cdot)\right\|^p _{L^p _{\gamma/p}}\,ds < +\infty
 \end{equation*}
and thus
 \begin{equation*}
 \forall t_0 >0,\ \  \exists \ t_1 \in (0,t_0) ; \quad
 \left\|f(t_1,\cdot)\right\|^p _{L^p _{\gamma/p}} < +\infty.
 \end{equation*}
Besides, the estimate~\eqref{eq:apriori:eqdif:weight} for $\eta =
\gamma/p$ gives the propagation of the $L^p _{\gamma/p}$-norm
starting from time $t_1 >0$. Since for $t\geq t_1$, the $L^1_s$
norms of $f$ are uniformly bounded, the argument can be iterated
to prove by induction (integrating in time the weighted
inequality~\eqref{eq:apriori:eqdif:weight}) that
 \begin{equation*}
 \forall \eta \ge 0, \ \forall t>0, \ \
 \left\|f(t,\cdot)\right\|_{L^p _\eta} <+\infty.
 \end{equation*}
The above argument is slightly formal since we worked with
quantities which are not a priori finite. It can however be made
rigorous and quantitative in the same manner as
in~\cite{Wenn:momt:97}.
\end{proof}

\Remark One could also prove the property of moment generation in $L^p$ directly,
without induction, by using the idea of Wennberg~\cite{Wenn:momt:97} of comparison
to a Bernoulli differential equation. Using the same estimates on $Q^+ _R$ and 
$Q^-$ as in~\eqref{eq:apriori:eqdif:weight}, convolution-like inequality~\eqref{eq:convol} on 
$Q^+ _S$, and H\"older inequality, one gets the following
 \begin{equation*}
 \frac{d\left\|f\right\|^p _{L^p _\eta}}{dt} \le
 C_+ \left\|f\right\|_{L^p _\eta} ^{p} -
 \frac{K_-}{C_p(f_0)} \left\|f\right\|_{L^p _{\eta}} ^{p(1+\lambda)}
 \end{equation*}
where $\lambda = \frac{\gamma}{\eta}$ and $C_p(f_0)$ stands for the uniform bound on
the $L^p$ norm of the solution. It gives an explicit bound on the $L^p$ moments of
the form
 \[ \forall t>0, \ \ \left\|f(t,\cdot)\right\|_{L^p _\eta}
 \leq \left[ \frac{A}{B\left(1-e^{-A \lambda t} \right)} \right]^{- \frac{\eta}{\gamma} } \]
where $A,B$ depend on $C_p(f_0)$ and an upper bound on $L^1$ moment of the solution
of high enough order. Notice that these bounds are not optimal: for example,
$\left\|f\right\|_{L^p _{\gamma/p}}$ has to be integrable as a function
of $t$, as $t \to 0^+$, as can be seen from our a priori differential inequality.

%%%%%%%%%%%%%%%%%%%%%%%%%%%%%%%%%%%%%%%%%%%%%%%%%%%%
%%%%%%%%%%%%%%%%%%%%%%%%%%%%%%%%%%%%%%%%%%%%%%%%%%%%
\subsection{Generalization: propagation of $H^k$ estimates for $k \in \N$}

Here we follow the same strategy on the differentiated equation in order to get
uniform bounds in Sobolev spaces $H^k$ for $k \in \N$. This method
seems to fail for spaces $H^k$ with $k$ non-integer, 
%seems to : mieux vaut rester prudent !
because fractional derivatives do not behave ``bilinearly'' with
respect to the collision operator. 
Moreover we state our results only for ``power law'' kinetic collision 
kernels. This restriction is made for convenience, and can probably be
relaxed at the price of some more work.

 \begin{theorem}\label{theo:prop:sobol}
 Let $B(v-v_*,\sigma) = |v-v_*|^\gamma \, b(\cos \theta)$ ($\gamma \in (0,2)$) satisfy
 assumptions~\eqref{eq:hypcut-off}, \eqref{eq:hypang2}, 
 \eqref{eq:hyplbB} and~\eqref{eq:hypL^2}, let $\eta \in \R$,
 and let $f_0$ be a nonnegative function in
 $L^1_2$. Then the unique solution $f$ of the
 Boltzmann equation with initial datum $f_0$ satisfies, for any
 multi-index $\nu$, the estimate
  \begin{equation*} \label{eq:diffineq:sobol}
  \frac{d}{dt} \left\|\partial^\nu f\right\|^2 _{L^2 _\eta} \le
  C_+ \left\|\partial^\nu f\right\|_{L^2 _\eta}  -
  K_- \left\|\partial^\nu f\right\|_{L^2 _{\eta + \gamma/2}} ^2
  \end{equation*}
 for some constants $C_+,K_->0$, which depend
 on $p$, $N$, $B$, on upper bounds on $\|f_0\|_{L^1_2} + H(f_0)$, on a
 lower bound on $\|f_0\|_{L^1}$ and
 on $L^2 _{\eta + 1 +\gamma}$ norms on derivatives of $f$ of order
 strictly less than $|\nu|$.

 In particular for any $k \in \N$, there is an explicit constant $C_k (f_0)$, only
 depending on $B$, on an upper bound on
 $\|f_0\|_{L^1_2}+\|f_0\|_{H^k _{k(1+\gamma)}}$,
 and on a lower bound on $\|f_0\|_{L^1}$, such that
  \[ \forall t\geq 0, \qquad
  \|f(t,\cdot)\|_{H^k} \leq C_k (f_0). \]
 Moreover, for any $t >0$ and any $\kappa > 0$, we know that
 $f(t,\cdot) \in H^k_\kappa(\R^N)$. More precisely, for any $t_0>0$,
  \[ \sup_{t\geq t_0} \|f(t,\cdot)\|_{H^k _\kappa} < +\infty. \]
 This bound can be computed in terms of
 $B$, an upper bound on $\|f_0\|_{L^1_2}+\|f_0\|_{H^k}$,
 a lower bound on $\|f_0\|_{L^1}$, and a lower bound
 on $t_0$.
 \end{theorem}

\begin{proof}[Proof of Theorem~\ref{theo:prop:sobol}]
%In order to make the proof understandable, we shall first do it in
%the simplified framework $k=1$. At the end we shall
%say briefly how to deal with the general case. 
Again we only prove the a priori differential inequality:
let us consider a given partial derivative $\partial^\nu f$ of $f$.
 \begin{align*}
 \frac{1}{2}\,
 \frac{d\left\|\partial^\nu f \right\|^2 _{L^2 _\eta}}{dt} &=
 \int \partial^\nu f \, \partial^\nu Q^+(f,f) \, \langle v \rangle^{2 \eta} \, dv 
 - \int \partial f \, \partial Q^-(f,f) \, \langle v \rangle^{2 \eta} \, dv \\
 &= \int \partial^\nu f \, \partial^\nu Q^+ \langle v \rangle^{2 \eta} \, dv 
 - \int (\partial^\nu f)^2 A \ast f \, \langle v \rangle^{2 \eta} \,dv \\
 &- \sum_{0 < \alpha \le \nu} \begin{pmatrix} \nu \\ \alpha \end{pmatrix}
 \int \partial^\nu f \, \partial^{\nu-\alpha} f \, \partial^\alpha(A \ast f) \, \langle v \rangle^{2 \eta} \, dv
 \end{align*}
where $A(z) = \|b\|_{L^1 (\ens{S}^N)} \Phi(z) = \cst |z|^\gamma$ here. 
% \begin{align*}
% \frac{1}{2}\,
% \frac{d\left\|\partial f \right\|^2 _{L^2}}{dt} &=
% \int \partial f \, \partial Q^+(f,f) \, dv - \int \partial f \, \partial Q^-(f,f)\,dv \\
% &= \int \partial f \, \partial Q^+ _S (f,f) \, dv +  \int \partial f \, \partial Q^+ _R(f,f) \, dv \\
% & - \int (\partial f)^2 A \ast f \,dv - \int \partial f \, f \, \partial(A \ast f) \,dv
% \end{align*}
For the first term we apply the regularity theorem~\ref{theo:reg:Q^+:BD} : since $(N-1)/2 \ge 1$, 
it implies
 \[ \left\|\partial^\nu Q^+ (f,f) \right\| _{L^2 _\eta}
 \le \Cbd \, \left[ \left\| f \right\|_{L^1 _{\eta + \gamma}} \left\| f \right\|_{L^1 _{\eta + \gamma}} 
                    + \left\| f \right\|_{H^{\nu'} _{\eta + \gamma+1}} 
                        \left\| f \right\|_{H^{\nu'} _{\eta + \gamma+1}} \right] \]
where $\nu'$ is a multi-index satisfying $|\nu'| < |\nu|$, and thus
 \[ \int \partial^\nu f \, \partial^\nu Q^+ (f,f) \, \langle v \rangle^{2 \eta} \, dv \le
 C_1 \, \left\| \partial^\nu f \right\|_{L^2 _\eta} \]
with $C_1$ depending on the $L^2 _{\eta + \gamma+1}$ norm on derivatives of $f$ of order
strictly lower than $\nu$ and the $L^1 _{\eta + \gamma}$ norm of $f$.

%For the second term the inequality used in the proof of 
%Theorem~\ref{theo:gain:reg} yields
% \[ \int \partial f \, \partial Q^+ _R(f,f)\,dv \le
% C_2 \, \left\| f \right\|_{L^1 _2} \left\| \partial f \right\|^2 _{L^2 _{\gamma/2}} \]
%with $C_2$ going polynomially to zero as $m,n$ go to infinity.
By Proposition~\ref{prop:minor}, the second term is bounded by
 \begin{equation} \label{eq:dampA}
 - \int (\partial f)^2 A \ast f \, \langle v \rangle^{2 \eta} \, dv
 \leq -K_0 \left\|\partial f\right\|^2 _{L^2 _{\eta + \gamma/2}}.
 \end{equation}
Finally for the third and last term, we split $A$ in $A_S + A_R$ where for $j \in \N$
 \[ A_S = \left( \widetilde{\Theta}_j \ast 1_{|v| \ge 2/j} \right) A \,, \qquad A_R = A - A_S \]
(notice that here we only need to isolate the singularity at zero
relative velocity). \\ 
For the smooth part,
 \[ \left\| \partial^\alpha(A_S \ast f) \right\|_{L^\infty} =
 \left\| (\partial^\alpha A_S) \ast f \right\|_{L^\infty} \le
 \left\| \partial^\alpha A_S \right\|_{L^\infty _{-(\gamma-1)^+}} \left\| f \right\|_{L^1 _{(\gamma-1)^+}} \]
($\left\| \partial^\alpha A_S \right\|_{L^\infty _{-(\gamma-1)^+}} < + \infty$ since $|\alpha| \ge 1$) 
and thus
 \[ \int \partial^\nu f \, \partial^{\nu-\alpha} f \, \partial^\alpha(A_S \ast f) 
 \, \langle v \rangle^{2 \eta} \, dv \le 
 C \, \left\| f \right\|_{L^1 _{(\gamma-1)^+}}
 \left\| \partial^\nu f \right\|_{L^2 _\eta} \left\| \partial^{\nu-\alpha}  f \right\|_{L^2 _\eta} \le 
 C_2 \, \left\| \partial^\nu f \right\|_{L^2 _\eta} \]
with $C_2$ depending on on the $L^2 _\eta$ norm on derivatives of $f$ of order
strictly less than $|\nu|$ and the $L^1 _{(\gamma-1)^+}$ norm of $f$. \\
For the remainder term,
 \[ \left\| \partial^\alpha (A_R \ast f) \right\|_{L^\infty} =
 \left\| A_R \ast (\partial^\alpha f) \right\|_{L^\infty} \le
 \left\| A_R \left\|_{L^2} \right\| \partial^\alpha f \right\|_{L^2} \]
and thus
 \[ \int \partial^\nu f \, \partial^{\nu-\alpha} f \, \partial^\alpha(A_R \ast f) \,dv \le 
 C_3 \, \left\| \partial^\nu f \right\| _{L^2} \]
if $\alpha < \nu$, with $C_3$ depending on the $L^2$ norm on derivatives of $f$ of order
strictly lower than $\nu$ and the $L^1$ norm of $f$, or 
 \[ \int \partial^\nu f \, \partial^{\nu-\alpha} f \, \partial^\alpha(A_R \ast f) \,dv \le 
 C_3 \, \left\| \partial^\nu f \right\|^2 _{L^2} \]
if $\alpha = \nu$, with $C_3$ depending on the $L^2$ norm of $f$. 
In the second case as $C_3$ goes to zero when $j$ goes to infinity, the term can be 
damped by the second one thanks to~\eqref{eq:dampA}. This shows that
 \[ \frac{1}{2}\,
 \frac{d\left\|\partial^\nu f \right\|^2 _{L^2}}{dt} \le C_+ \left\| \partial^\nu f \right\| _{L^2} -
 K_- \left\| \partial^\nu f \right\|^2 _{L^2 _{\gamma/2}} \]
and the proof is complete.

Then the proof of propagation of $H^k$ norm is made by induction. The proof of 
moments appearance is made first by propagating the $H^k _{-k(\gamma+1)}$ norm, then 
using interpolation with the $L^1$ moments.

\end{proof}

\Remark To get $W^{k,p}$ bounds when $p$ is different from $2$,
the strategy above could still apply, although with more
complications. The idea would be to prove an a priori differential inequality 
similar to~\eqref{eq:diffineqp} on each derivative. 
One should use the decomposition $Q^+ = Q^+ _S + Q^- _R$. 
To deal with the regular part one should now use
Corollary~\ref{coro:gain:leb} instead of Theorem~\ref{theo:reg:Q^+:BD} 
on each term of the Leibniz formula; and to deal with the 
remainder part one should use estimate~\eqref{eq:dual:sym:deriv},
together with the rough estimate
 \[ \| f \|_{L^1 _\eta} \leq C(\var) \| f \|_{L^p _{\eta + N/p' + \var}} \]
for $\var >0$. Moreover the weight exponent in the
assumptions become much higher.

%%%%%%%%%%%%%%%%%%%%%%%%%%%%%%%%%%%%%%%%%%%%%%%%%%%%
%%%%%%%%%%%%%%%%%%%%%%%%%%%%%%%%%%%%%%%%%%%%%%%%%%%%
%%%%%%%%%%%%%%%%%%%%%%%%%%%%%%%%%%%%%%%%%%%%%%%%%%%%

\section{Propagation of smoothness and singularity via Duhamel formula} \label{sec:propag}
\setcounter{equation}{0}

The aim of this section is to study the propagation of
smoothness and singularity for the solutions of the Boltzmann equation.
Throughout the section, we shall consider a given collision kernel
$B$, satisfying assumptions~\eqref{eq:hypcut-off},~\eqref{eq:hyprad},
~\eqref{eq:hypang2},~\eqref{eq:hyplbB},~\eqref{eq:hypL^2}.

%%%%%%%%%%%%%%%%%%%%%%%%%%%%%%%%%%%%%%%%%%%%%%%%%%%%
%%%%%%%%%%%%%%%%%%%%%%%%%%%%%%%%%%%%%%%%%%%%%%%%%%%%

 \subsection{Preliminary estimates}

From now on, explicit computations become rather long and we shall
try to be as synthetical as possible; so we will not keep track of
exact constants. However, all the proofs remain completely
explicit and there would be no conceptual difficulty in extracting
exact constants.

Our results in the sequel are based on two kinds of estimates.
First, a result of stability in $L^1$ for the solution of the
Boltzmann equation with cut-off and hard potential. Secondly, some
smoothness estimates on the Duhamel representation formula.
\medskip

\bul The stability result in $L^1$ which we use
is an immediate consequence of the
estimates in~\cite{Wenn:momt:94} and in~\cite{Gust:L^p:86}.
We do not search here for an optimal version. As in the sequel,
we shall use the shorthands $f_t = f(t,\cdot)$.

 \begin{lemma}\label{lem:stab}
 Let $f,g$ be two solutions of the Boltzmann equation
 belonging to $L^1 _{2+\gamma} \cap L\log L$.
 Then there exists a constant $C >0$,
 only depending on $B$, such that for all $0 \le k \le 2$ and $t \ge 0$
  \begin{equation*}\label{eq:stab:der:L^1}
  \frac{d}{dt}\left\|f_t-g_t\right\|_{L^1 _k} \le
  C \left\|f_t - g_t \right\|_{L^1 _k}
  \left\|f_t + g_t \right\|_{L^1 _{k +\gamma}}.
  \end{equation*}
 In particular, as $\left\|f_t + g_t \right\|_{L^1 _{k
 +\gamma}}$ is bounded uniformly with respect to $t$ thanks to the assumption, and we have
 the stability estimate
  \begin{equation*} \label{eq:stab:L^1}
  \left\|f_t - g_t\right\|_{L^1 _k} \le
  \left\|f_0 - g_0 \right\|_{L^1 _k} e^{\Cstab t},
  \end{equation*}
 where $\Cstab$ only depends on $B$, $\|f_0\|_{L^1_{k+\gamma}}$ and $\|g_0\|_{L^1_{k+\gamma}}$.
 \end{lemma}

\bul Next, we introduce the well-known Duhamel representation formula
for the Boltzmann equation,
 \begin{equation}\label{eq:form:duham}
 \forall \, t \ge 0, \ v \in \R^N, \ \
 f(t,v) = f_0 (v) e^{-\int_0 ^t Lf (s,v) \, ds }
 + \int_0 ^t Q^+ (f,f)(s,v)  e^{-\int_s ^t Lf (\tau,v) \, d\tau } \, ds
 \end{equation}
where $Lf$ stands for $A \ast f$ and $A(z) = \left\|b\right\|_{L^1
(\ens{S}^N)} \Phi(|z|)$. This formula is well-adapted to the study
of smoothness issues because it expresses the solution in terms of
the initial datum and the regularizing operators $Q^+$ and $L$.

For $s \le t$, we set
 \begin{equation*}
 F(s,t,v) = \int_s ^t Lf(\tau,v) \, d\tau, \qquad
 G(s,t,v) = e^{-F(s,t,v)}.
 \end{equation*}

We shall prove several estimates on these functions. We look for
uniform (with respect to time) estimates, which
leads us to allow a ``loss'' on the weight exponent.
 \begin{proposition}\label{prop:estim:duhamel}
 Let $\alpha, \beta>0$ be such that
 $A \in H^\alpha _{-\beta}$. Let $\alpha'=\min(\alpha,(N-1)/2)$, and let $\delta=\beta+\gamma+1$.
 Then, there is a constant $\Cduh$ such that for all $k, \eta\geq 0$,
  \begin{equation} \label{eq:reg1}
  \left\|\int_0 ^t Q^+(f,f)(s,\cdot) \, G(s,t,\cdot) \, ds
  \right\|_{H^{k+\alpha'}_\eta}  \le \Cduh \sup_{0 \le \tau \le t}
  \left\| f(\tau,\cdot) \right
  \|_{H^k _{\eta+\delta}}^{\lceil k+\alpha \rceil + 2},
  \end{equation}
 and
  \begin{equation} \label{eq:reg2}
  \left\|f_0(\cdot) \, G(0,t,\cdot) \right\|_{H^{k} _\eta}
  \le \Cduh \, e^{-K' t}
  \left\|f_0(\cdot) \right\|_{H^{k} _{\eta+\beta}}
  \sup_{0 \le \tau \le t}
  \left\| f(\tau,\cdot) \right\|^{\lceil k \rceil}
  _{H^{k-\alpha'} _\beta},
  \end{equation}
 with $0< K' < K$ where $K>0$ is the constant in~\eqref{eq:lowbdQ^-}.
 \end{proposition}

\Remark Under our general assumptions, a possible choice
of $\alpha,\beta$ is $\alpha=\gamma$, $\beta = N/2+\gamma+\var$,
$\var>0$. For $\Phi(|z|)=|z|^\gamma$, it would be possible to
take $\alpha=\gamma+N/2-\var$, for any $\var>0$.
\smallskip

\begin{proof}[Proof of Proposition~\ref{prop:estim:duhamel}]
We start with some preliminary estimates on $L$, $F$ and $G$.
As a consequence of Cauchy-Schwarz inequality,
we find that for all $k\geq 0$,
 \begin{equation*} \label{eq:estim:Lf}
 \left\|Lf\right\|_{W_{-\beta} ^{k+\alpha,\infty}} \le
 C_1 \| f \|_{H^k _\beta}.
 \end{equation*}
It follows that
 \[  \left\|F(s,t,\cdot)\right\|_{W_{-\beta} ^{k+\alpha,\infty}} \le
 C_1 \, \sqrt{t-s} \,
 \left( \int_s ^t \| f(\tau,\cdot) \|^2 _{H^k _\beta} \, d\tau \right)^{1/2}. \]
Combining this with the estimate~\eqref{eq:lowbdQ^-:faible}, in
the form $Lf \geq K$, we deduce that
 \begin{eqnarray} \label{eq:estimG}  \nonumber
 \left\| G(s,t,\cdot)\right\|_{W_{-\beta} ^{k+\alpha,\infty}} &\le& C_1 \,
 \sqrt{t-s} \, e^{-K(t-s)} \left(\int_s ^t \| f(\tau,\cdot) \|^2 _{H^k _\beta} \, d\tau 
 \right)^{\frac{\lceil k+\alpha \rceil}{2}} \\
 &\le& C_2 \, e^{-K' (t-s)} \, \sup_{s \le \tau \le t}
 \left\| f(\tau,\cdot) \right\|^{\lceil k +\alpha \rceil}  _{H^k _\beta}
 \end{eqnarray}
with $0< K ' < K$. 

Now we use the following simple lemma to exchange a time integral and a
$H^S _P (\R^N_v)$ norm:
 \begin{lemma}\label{lem:estim:intervert}
 Let $Z(s,v)$ be a function on $\R_+ \times \R^N$ and $S,P \in R$, then 
 for any $\lambda > 0$
  \begin{equation*}\label{eq:estim:intervert}
  \left\| \int_0 ^t Z(s,\cdot) \, ds \right\|_{H^S _P}
  \le \frac{1}{\sqrt{\lambda}} \left( \int_0 ^t e^{\lambda(t-s)}
  \left\| Z(s,\cdot) \right\|^2 _{H^S _P} \, ds \right)^{1/2}.
  \end{equation*}
 \end{lemma}
This lemma is an immediate consequence of
the Cauchy-Schwarz inequality with the weight
$e^{\frac{\lambda}{2}(t-s)}$, after passing to Fourier variables.
The choice of the exponential function is arbitrary; we used it
because it is convenient for the sequel.

As a consequence, we have (recall that $\alpha'=\min(\alpha,(N-1)/2)$)
 \[ \left \| \int_0^t Q^+(f_s,f_s)\, G(s,t,\cdot)\,ds
 \right \|_{H^{k+\alpha'}_\eta} \]
 \[ \leq C \left ( \int_0^t e^{K'(t-s)} \Bigl \|Q^+(f_s,f_s)
 \,G(s,t,\cdot) \Bigr \|^2_{H_\eta^{k+\alpha'}} \,ds \right )^{1/2} \]
 \[\leq C \left ( \int_0^t e^{K'(t-s)}
 \|Q^+(f_s,f_s)\|^2_{H^{k+\alpha'}_{\eta+\beta}}
 \|G(s,t,\cdot) \|^2_{W^{k+\alpha',\infty}_{-\beta}}\,ds \right )^{1/2}. \]
% \[ \leq C \left ( \int_0^t e^{K(t-s)}
% \|Q^+(f_s,f_s)\|^2_{H^{k+\alpha'}_{\eta+\beta}}
% \|G(s,t,\cdot)\|^2_{W^{k+\alpha',\infty}_{-\beta}}\,ds \right )^{1/2}. \]
At this stage we apply Theorem~\ref{theo:reg:Q^+:BD} and
estimate~\eqref{eq:estimG}, to get a bound like
 \[ C \left [ \int_0^t e^{K'(t-s)} \, 
 \|f_s\|^4_{H^k_{\eta+\beta+\gamma+1}}
 e^{-2K'(t-s)} \, \left( \sup_{s \le \tau \le t}
 \left\| f(\tau,\cdot) \right\|^{\lceil k +\alpha \rceil}  _{H^k _\beta}\right)^2 \, ds \right]^{1/2} \]
% \left ( \int_s^t \|f_\tau\|^2 _{H^{k}_{\beta}} \right )^{
% \lceil k+\alpha \rceil}\,ds
% \right ]^{1/2}\]
 \[\leq C \left ( \int_0^t e^{-K'(t-s)} ds \right )^{1/2}
 \sup_{0\leq s\leq t} \|f(s,\cdot)\|^{\lceil k+\alpha\rceil+2} _{H^k_{\eta+\beta+\gamma+1}}
 \leq C \sup_{0\leq s\leq t} \|f(s,\cdot)\|^{\lceil k+\alpha\rceil+2} _{H^k_{\eta+\beta+\gamma+1}}. \]
This concludes the proof of~\eqref{eq:reg1}.

The proof of~\eqref{eq:reg2} is performed in a similar way, using estimate~\eqref{eq:estimG} with $s=0$. 
\end{proof}

%%%%%%%%%%%%%%%%%%%%%%%%%%%%%%%%%%%%%%%%%%%%%%%%%%%%
%%%%%%%%%%%%%%%%%%%%%%%%%%%%%%%%%%%%%%%%%%%%%%%%%%%%

 \subsection{Propagation of regularity}

As soon as we have uniform bounds on $L^2$ moments, the
Duhamel representation~\eqref{eq:form:duham} together with
Proposition~\ref{prop:estim:duhamel} imply some
uniform bounds on $f$ in Sobolev spaces, provided that the initial
datum itself belong to such a space. With respect to the method used 
for proving Theorem~\ref{theo:prop:sobol}, the improvement here
is that we are able to treat $H^s$ regularity for any $s \in \R_+$. 
Here is a precise theorem, definitely not optimal.

 \begin{theorem}\label{theo:prop:reg}
 Let $0 \le f_0\in L^1_2$ be an initial datum
 with finite mass and kinetic energy, and let $f$ be the unique
 solution preserving energy. Then for all $s>0$ and $\eta \ge \beta$,
 there exists $w(s)>0$ (explicitly $w(s)=\delta \lceil s/\alpha' \rceil$) such that
  \[ f_0\in H^s_{\eta+w} \Longrightarrow
  \sup_{t\geq 0} \|f(t,\cdot)\|_{H^s _\eta} < +\infty.\]
 \end{theorem}

\Remark This theorem is not so strong as the decomposition theorem
below, because of the strong moment assumption. It is quite likely
that the restriction about $w$ could be relaxed with some more work.
A sufficient condition for this moment
assumption to be automatically satisfied, is that all the $L^1$
moments of $f_0$ be finite. Of course we know that for $t\geq
t_0$, this is always the case; but this is a priori not sufficient
to conclude. Nevertheless it gives by interpolation the following
result: under the same assumptions, as soon as $f_0 \in H^s
_\eta$, $f_t$ belongs to $H^s _\eta$ for any $t>0$. The constant is
explicit, is uniformly bounded for $t>t_0$ for any $t_0$, and
blows-up like an inverse power law of $t_0$ as $t_0\to 0^+$.

\begin{proof}[Proof of Theorem~\ref{theo:prop:reg}]
Let $n \in \N$ be such that $n \alpha' \geq s$ ($n=\lceil
s/\alpha' \rceil$). Let $w(s)=\delta \lceil s/\alpha' \rceil$. The
proof is made by an induction comprising $n$ steps, 
proving successively that $f$ is uniformly bounded in 
$H^{i \alpha'} _{\eta + \frac{n -i}{n} w}$
for $i= 0,1,...,n$. The above-mentioned argument is used in each step.

Let us write the induction. The initialisation for $i=0$, i.e $f$
uniformly bounded in $L^2 _{\eta + w}$ is proved by
Theorem~\ref{theo:prop:L^p} and the more general
equation~\eqref{eq:apriori:eqdif:weight}. Now let $0 < i \le n$
and suppose the assumption is satisfied for all $0 \le j < i$.
Then proposition~\ref{prop:estim:duhamel} implies
 \[
 \|f_0(\cdot) \, G(0,t,\cdot) \|_{H^{i \alpha'} _{\eta + \frac{n -i}{n} w}}  
 \le C_2 \, e^{-K't}
 \left\|f_0(\cdot) \right\|_{H^{i \alpha'} _{\eta + \frac{n - i}{n} w +
 \beta}}
 \> \sup_{ 0 \le \tau \le t}
 \left\| f(\tau,\cdot) \right\|_{H^{(i-1)\alpha'} _\beta} ^{\lceil i  \alpha' \rceil}.
 \]
We know from the previous subsection that
 \[
 \left \|\int_0 ^t Q^+(f,f)(s,\cdot) \,
 G(s,t,\cdot) \, ds \right\|_{H^{i \alpha'} _{\eta + \frac{n -i}{n}w}} 
 \le  C_3 \> \sup_{0 \le \tau \le t}
 \left\| f(\tau,\cdot) \right\|_{H^{(i-1)\alpha'}
 _{\eta + \frac{n -i}{n}w + \delta}} ^{\lceil i \alpha' \rceil +2}
 \]
Moreover as $\beta \le \delta \le w/n$ and $i \ge 1$,
 \begin{eqnarray*}
 \beta &\le& \eta + \frac{n-(i-1)}{n}w\\ 
 \eta + \frac{n-i}{n}w + \delta &\le& \eta + \frac{n -(i-1)}{n} w
 \end{eqnarray*}
and thus, using the induction assumption for $i-1$, $f$ is uniformly bounded in $H^{k \alpha'} _{\eta
+\frac{n-i}{n}w}$ and the proof is complete.
\end{proof}
%%%%%%%%%%%%%%%%%%%%%%%%%%%%%%%%%%%%%%%%%%%%%%%%%%%%
%%%%%%%%%%%%%%%%%%%%%%%%%%%%%%%%%%%%%%%%%%%%%%%%%%%%

 \subsection{The decomposition theorem}

Here we shall give a precise meaning to the idea that the
Boltzmann equation with cut-off propagates both smoothness and
singularities, but makes the amplitude of the singular part go to
zero as time $t$ go to infinity. To this purpose, we shall look for some
iterated versions of the Duhamel representation~\eqref{eq:form:duham}. 

 \begin{theorem}\label{theo:dec}
 Let $0 \le f_0 \in L^1 _2 \cap L^2$ and $f$ be the unique energy-preserving
 solution of the Botzmann equation with initial datum $f_0$, and let
 $s\geq 0$, $q\geq 0$ be arbitrarily large.
 Let $\tau>0$ be arbitrarily small. Then, for all $t\geq \tau$,
 $f$ can be written $f^S+f^R$, where $f^S$ is nonnegative,
 and
  \begin{equation*}\label{eq:decomp:controle}
  \left\{
   \begin{array}{l}
   \displaystyle \sup _{t \ge \tau} \left\|f_t ^S \right\|_{H^s_q \cap L^1 _2} < +\infty \\ \\
   \displaystyle \forall t \ge \tau,\> \forall k>0,\> \exists
   \lambda=\lambda(k)>0; \ \ \left\|f_t ^R \right\|_{L^1 _k}
   = O\left(e^{-\lambda t}\right) .
   \end{array}
  \right.
  \end{equation*}
 All the constants in this theorem can be computed in terms of
 the mass, energy and $L^2$-norm of $f_0$, and $\tau$.
 \end{theorem}
\Remark The idea of such a decomposition is reminiscent of Wild sums in the case 
of Maxwellian molecules. Also partial results in this direction 
were obtained in Wennberg~\cite{Wenn:rado:94} and Abrahamsson~\cite{Abra:99}. In these 
cases the gain of regularity in the second term of the Duhamel formula was iterated just once (or twice in~\cite{Abra:99} for a gain of integrability), 
and thus the regularity was limited to $H^{(N-1)/2}$ 
essentially. For hard potentials the obstacle to iterate the Duhamel formula as in the 
Maxwellian case is the strong non-linearity of the 
decomposition. Here we bypass this difficulty by the strategy of starting new flows at each step of the iteration. 
\medskip
 
\begin{proof}[Proof of Theorem~\ref{theo:dec}]
We first note that moment estimates imply bounds in $L^1_k$ for all $k\geq 0$,
and therefore the only problems are the gain of regularity for the
smooth part and the exponential decrease for the remainder part.

The idea of the proof is a use of the Duhamel formula to decompose the flow associated
with the equation into two parts, one of which is more regular than the initial
datum, while the amplitude of the other decreases exponentially fast with time. We
shall use this repeatedly to progressively increase the smoothness: after a while,
we start again a new flow having the smooth part of the previous solution as initial
datum. And so on. Of course, each time we start a new flow, we shall depart from the
true solution, since the initial datum is not the real solution. However, we can use
the stability theorem (Lemma~\ref{lem:stab}) to control the error.

The times at which we start the new flows are chosen in such a way
that the decay of the non-smooth part (measured by the constant
$\Cdec$) balances the divergence of the solutions (measured by the
constant $\Cstab$).
The idea is summarized in figure~\ref{fig1}. Each node of the
tree corresponds to a time where we start a new solution of the
Boltzmann equation taking for initial data the ``smooth part'' of
the previous solution.
In the aim to achieve the goal of balancing the effect of the
divergence of the solutions thanks to the exponential decaying of
the first term in the Duhamel formula, it is necessary that the
decomposition tree ends precisely at the time $t$ we are looking
for a decomposition of the solution. Note that for different
$t$, the functions $f^S _t$ constructed below {\em do not belong
to the same flow}.
\begin{figure}
\epsfysize=8cm
$$\epsfbox{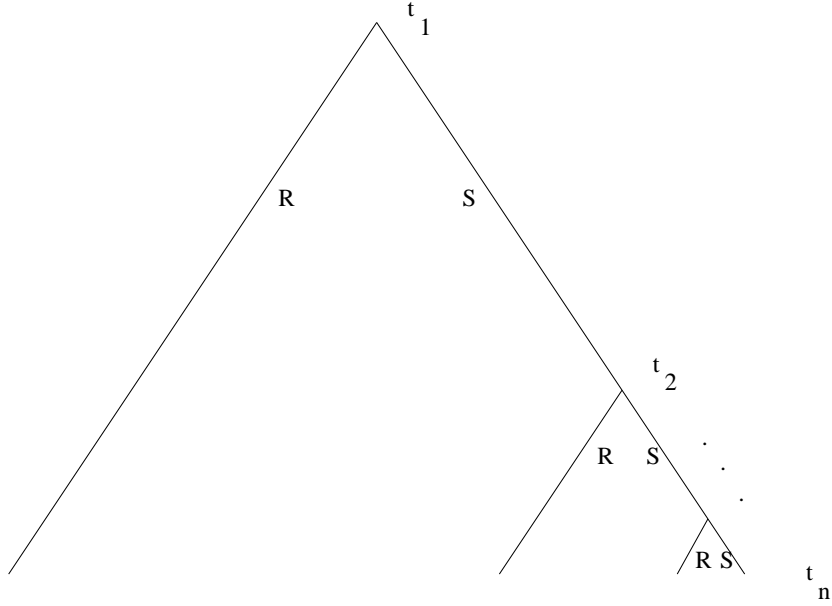}$$
\caption{Decomposition of the solution}\label{fig1}
\end{figure}

Let us implement this idea more precisely.
By Theorem~\ref{theo:prop:L^p}, we have a uniform $L^2$ bound on
the solution $f$, and for a given $t_0 > 0$, we also know that
all the $L^2$-moments are uniformly bounded (see subsection~\ref{sec:Lp}).
Let $n\geq 1$, to be thought of as the number of times we wish to
apply the semigroup; we choose $n$ in such a way that $n\alpha' >
k$, where $k$ is the degree of smoothness which we are looking
for, and $\alpha'$ is the degree of regularization appearing in
Proposition~\ref{prop:estim:duhamel}. Let $\tau'\in (0,\tau)$ be
arbitrary, say $\tau'=\tau/2$. Let us set $t_n=t \geq \tau$,
$t_{-1}=\tau'$, and define inductively
(forwards) $t_i$ for $0\leq i\leq n-1$ by
 \begin{equation*} \label{eq:def:tps:dec}
 t_i = t_{i-1} + \mu ( t_n - t_{i-1})
 \end{equation*}
where $\mu \in (0,1)$ satisfies
 \begin{equation} \label{eq:cond:mu}
 \mu > \frac{\Cstab}{\Cstab + K'}
 \end{equation}
($K'$ is the constant of exponential decrease in~\eqref{eq:reg2}). 
Let us denote $f_0,f_1,...,f_n$ the solutions constructed as
explained above: $f_0 = f$ is the solution we are studying, $f_1$
is the solution for $t \ge t_0$ of the Boltzmann equation starting
from the ``initial datum''
 \begin{equation*}
 \int_0 ^{t_0} Q^+(f_0,f_0)(s,\cdot) \, e^{-\int_s ^{t_0} Lf_0 } \, ds
 \end{equation*}
at time $t_0$, $f_2$ is the solution for $t \ge t_1$
of the Boltzmann equation starting from
 \begin{equation*}
 \int_0 ^{t_1 - t_0} Q^+(f_1,f_1)(s,\cdot) \,
 e^{-\int_s ^{t_1 -t_0} Lf_1 } \, ds
 \end{equation*}
at time $t_1$, etc. More generally, for $2 \le i \le n-1$,
$f_{i+1}$ is the solution for $t \ge t_i$ of the Boltzmann equation
starting at
 \begin{equation*}
 \int_0 ^{t_i - t_{i-1}} Q^+(f_i,f_i)(s,\cdot) \,
 e^{-\int_s ^{t_i -t_{i-1}} Lf_i } \, ds
 \end{equation*}
at time $t_i$. Of course this sequence is well-defined,
since at each node, the ``smooth'' part of the solution that we
take as a new initial data is nonnegative, lies in $L^1_2 \cap L^2$
and has all its $L^2$-moments bounded.

The $n$-times iteration of estimate~\eqref{eq:reg1} together with
the theorem of propagation of regularity~\ref{theo:prop:reg}
easily implies a bound on the $H^{n\alpha'}$ norm on $f_n
(t_n,\cdot)$ which is uniform in $t_n \geq \tau$, and only depends
on $\tau'$, $n$, and on the mass, energy and $L^2$-moments of $f$.
So let us set
 \begin{equation*}
 f^S (t_n, \cdot) \equiv f_n(t_n,\cdot)
 \end{equation*}
and
 \begin{equation*}
 f^R (t_n, \cdot) \equiv f(t_n,\cdot) - f^S (t_n,\cdot).
 \end{equation*}
This construction can be made for all $t_n \ge \tau$; thus our decomposition
is well-defined for all $t \geq \tau$. It remains to prove
that $f^R$ is exponentially decaying as $t\to\infty$. For this we
write
 \begin{eqnarray*}
 \left\|f^R _t\right\|_{L^1} &=&
 \left\|f^R _{t_n} \right\|_{L^1} \\ 
 &\le&
 \sum_{i=0} ^{n-1}
 \left\|f^{i+1}  _{t_n} - f^i _{t_n}\right\|_{L^1} \\ 
 &\le& \sum_{i=1} ^n e^{\Cstab (t_n-t_i)}
 \left\|f^{i+1} _{t_i} - f^i _{t_i}\right\| \\ 
 &\le& C \sum_{i=1} ^n e^{\Cstab (t_n-t_i)}
 e^{-K'(t_i-t_{i-1})} \\ \nonumber
 &\le& C \sum_{i=1} ^n e^{(1-\mu)^{i-1} (t_n-t_0)
 \left(\Cstab (1-\mu) -K' \mu \right)} \\ 
 &\le& C \sum_{i=1} ^n e^{(1-\mu)^{n-1} (t_n - \tau')
 \left(\Cstab (1-\mu) -K' \mu \right)} \\ 
 &\le& C n e^{ - (t_n-\tau') (1-\mu)^{n-1}
 \left(K' \mu - \Cstab (1-\mu) \right)}
 \end{eqnarray*}
which gives the result: if one set
 \begin{equation*}
 0 < \Cdec < (1-\mu)^{n-1} \left(K' \mu - \Cstab (1-\mu) \right)
 \end{equation*}
which is possible thanks to~\eqref{eq:cond:mu}, we have
 \begin{equation*}
 \left\|f^R _t\right\|_{L^1} \le C e^{-\Cdec t}.
 \end{equation*}
On the other hand, $f^R_t$ has all its $L^1_k$ norms bounded, for all $k$.
%La justification serait : des t_1 f possede tous les mts L^1
%qui sont donc propages par f et par les f_S construits a partir de
%donnees initiales plus petites que f. Donc par diffce, f_R
%a tous ses mts L^1 bornes. Je ne pensais pas le mettre, dis-moi
%ce que tu en penses.
%
%% OK
By elementary interpolation, it follows that all these $L^1_k$ norms are
decaying exponentially fast (the same holds true for all $L^p_k$
norms, whenever $p<2$).
\end{proof}

%\Remark
%As a consequence of the preceding theorem, when $f_0\in L^2$ one can give a
%simplified proof of the Pulvirenti-Wennberg result of appearance
%of a Maxwellian lower bound~\cite{PulvWenn:binf:97}. Indeed, the
%regularization property easily implies that after a short time
%$f_t$ is bounded below by $f_t^S$ (since  $f_t^R\geq 0$), which
%itself is bounded below by a multiple of the characteristic function
%of a ball. With this one can bypass the very technical lemma~3.1
%in~\cite[section~3]{PulvWenn:binf:97}.

%%%%%%%%%%%%%%%%%%%%%%%%%%%%%%%%%%%%%%%%%%%%%%%%%%%%
%%%%%%%%%%%%%%%%%%%%%%%%%%%%%%%%%%%%%%%%%%%%%%%%%%%%
%%%%%%%%%%%%%%%%%%%%%%%%%%%%%%%%%%%%%%%%%%%%%%%%%%%%

\section{Application to a problem of long-time behavior} \label{sec:longtime}
\setcounter{equation}{0}

Let us now show an application of Theorem~\ref{theo:dec}.
Here we shall extend a result proven for very smooth solutions, 
into a result which applies without smoothness assumption.

We start with the following statement, which is an immediate corollary 
of the main results in~\cite{Vill:cveq:TA}.

 \begin{theorem}\label{theo:cv:2001}
 Let $B$ satisfy assumptions~\eqref{eq:hyprad}, \eqref{eq:hyplbB}
 and~\eqref{eq:hypL^2}, together with the stronger lower bound assumption
  \begin{equation} \label{eq:hyplbB2}
  b(\cos\theta) \geq b_0 >0.
  \end{equation}
 Let $f_0$ be a nonnegative function in $L^1_2(\R^N)$. Without loss
 of generality, assume that
 $\int f_0 =1$, $\int f_0(v)v\,dv =0$, $\int f_0(v)|v|^2\,dv = N$,
 and denote by
  \[ M(v) = \frac{e^{-|v|^2}}{(2\pi)^{N/2}} \]
 the associated Maxwellian equilibrium.
 Let $f$ be an energy-preserving solution of the Boltzmann equation
 with initial datum $f_0$, satisfying
  \begin{equation} \label{Hks}
  \forall s, \ \forall k, \qquad C_{s,k} \equiv
  \sup_{t\geq t_0} \|f_t\|_{H^s _k} < +\infty,
  \end{equation}
 and
  \[ \forall t\geq t_0, \qquad f(t,v) \geq K_0 e^{-A_0|v|^{q_0}}. \]
 for some time $t_0>0$ and some positive constants $K_0, A_0, q_0$.
 Then, $\|f_t-M\|_{L^1} = O(t^{-\infty})$, in the sense that
 for all $\var>0$ there exists a constant $C_\var$, explicitly
 computable in terms of the above constants, and depending on $f$
 only via $t_0$, $K_0$, $A_0$, $q_0$ and
 and upper bound on $C_{k,s}$ for $k$ and $s$ large enough, such that
  \begin{equation} \label{ccltinfty}
  \|f_t- M \|_{L^1} \leq  C_\var t^{-1/\var}.
  \end{equation}
 \end{theorem}

\Remarks 1. Assumption~\eqref{eq:hyplbB2} is satisfied by the hard
spheres kernel, and can be considered as satisfactory for hard
potentials with cut-off (since they are satisfied by non-cutoff
potentials). Kernels like $|v-v_*|^\gamma$ ($0<\gamma<2$) satisfy
all the above assumptions.
\smallskip

2. Note that Theorem~\ref{theo:cv:2001} and Theorem~\ref{theo:cvg:eq} 
in the sequel are quantitative, which explains their interest 
even if exponential convergence to equilibrium has been proven by 
non-constructive approaches: see Arkeryd~\cite{Arke:stab:88} 
for the proof in the $L^1$ setting, and Wennberg~\cite{Wenn:stab:93} 
for the extension to the $L^p$ setting.  
\medskip

It is equivalent to require~\eqref{Hks} or to require uniform
bounds in all $H^k$ norms and in all $L^1_s$ norms. Therefore,
we see that known results of appearance of moments and Maxwellian
lower bound for hard potentials cover all the assumptions needed
for this theorem, except the $H^k$ bounds. If we apply
the propagation result of Theorem~\ref{theo:prop:reg}, we conclude that the
conclusion~\eqref{ccltinfty} holds true as soon as the initial
datum $f_0$ lies in all weighted Sobolev spaces. However, the
decomposition theorem of the previous paragraph will lead us to
a much stronger conclusion.

 \begin{theorem}\label{theo:cvg:eq}
 Let $f_0$ satisfy the same assumptions as in
 Theorem~\ref{theo:cv:2001} and $B$ satisfy~\eqref{eq:hypcut-off},
 ~\eqref{eq:hyprad},~\eqref{eq:hypang2},~\eqref{eq:hyplbB},~\eqref{eq:hypL^2}
 together with~\eqref{eq:hyplbB2}. Further assume that $f_0\in L^2(\R^N)$.
 Then the conclusion of Theorem~\ref{theo:cv:2001} holds true:
  \begin{equation*}
  \left\|f_t - M \right\|_{L^1} = O ( t^{-\infty}).
  \end{equation*}
 \end{theorem}

\begin{proof}[Proof of Theorem~\ref{theo:cvg:eq}]
First of all, let us pick a $t_0 > 0$. We know that the solution
$f_t$ satisfies a Maxwellian lower bound and moment estimates,
uniformly as $t\geq t_0$.

Let $\var>0$ be arbitrary, and let $k, s$ be such that $C_\var$
in Theorem~\ref{theo:cv:2001} only depends on a uniform upper bound on
$\|f\|_{H^k_s}$.
Let us make the decomposition
of Theorem~\ref{theo:dec} with $\tau=1$ and $s$. Then we know that
 \begin{equation*}
 f_t = f^S _t + f^R _t, \ \ \ \forall \ t \ge 1.
 \end{equation*}
Let $t_1 \geq 1$ be an intermediate time, to be chosen later. Let
us introduce $\tilde{f}_t$ the solution of the Boltzmann equation
starting from $f^S _{t_1}$ at $t=t_1$, and $\tilde{M}$ the
Maxwellian distribution associated with $\tilde{f}_{t_1} = f^S
_{t_1}$. Since $f^S _{t_1}$ is bounded in $H^k \cap L^1 _s$ 
by theorem~\ref{theo:dec}, $\tilde{f}_t$ is uniformly bounded in $H^k \cap L^1
_s$, and has a Maxwellian lower bound
for $t \ge t_2 > t_1$ where $t_2$ can be chosen arbitrarily (so
let us say $t_2 = t_1 +1$). After rescaling space (to reduce to
the case where $\tilde{f}$ has unit mass, zero average velocity
and unit temperature), Theorem~\ref{theo:cv:2001} implies
 \begin{equation*}
 \|\tilde{f} - \tilde{M} \|_{L^1}
 = O ( (t-t_1) ^{-\frac{1}{\var}})
 \end{equation*}
with explicit constants which do not depend on $t_1$ (they only
depend on the $\tau$ in the decomposition).

Now, thanks to the properties of the decomposition
 \begin{eqnarray*}
 \left\| f^R _t \right\|_{L^1 _\gamma} &=& O\left(e^{-\Cdec(t-\tau)}\right) \\ 
 &=& O\left(e^{-\Cdec(t-1)}\right) = O\left(e^{-\Cdec t}\right).
 \end{eqnarray*}
Moreover,
 \begin{eqnarray*}
 \left\|M-\tilde{M} \right\|_{L^1} &=& O\left(e^{-\Cdec(t_1-\tau)}\right) \\ 
 &=& O\left(e^{-\Cdec(t_1-1)}\right) = O\left(e^{-\Cdec t_1}\right).
 \end{eqnarray*}
Indeed, a simple computation shows that $\|\tilde{M}-M\|_{L^1}$
can be bounded in terms of $\|f-\tilde{f}\|_{L^1_2}$, which in turn
can be estimated in terms of $\|f^R_{t_1}\|_{L^1_2}$.

Next, the stability lemma~\ref{lem:stab} implies
 \begin{equation*}
 \left\|f_t - \tilde{f}_t \right\| \le
 C e^{ \Cstab (t-t_1)} \left\|f_{t_1} - \tilde{f}_{t_1} \right\|_{L^1 _\gamma}.
 \end{equation*} 
On the whole, we find
 \begin{eqnarray*} 
 \left\|f_t - M \right\|_{L^1} &\le&
 \left\|f_t - \tilde{f}_t \right\|_{L^1} +
 \left\|\tilde{f}_t - \tilde{M} \right\|_{L^1} +
 \left\|\tilde{M} - M \right\|_{L^1} \\ 
 &\le& C \left( e^{\Cstab(t-t_1)} e^{-\Cdec t_1} +
 (t-t_1)^{- \frac{1}{\var} } + e^{-\Cdec t_1} \right)
 \end{eqnarray*}
It remains to choose $t_1 \gg (t-t_1)$ in order to compensate for the exponential
divergence allowed by the stability lemma~\ref{lem:stab}. More precisely, if
$\Cstab(t-t_1) = \frac{\Cdec}{2} \, t_1$ (i.e $t_1 = \Cstab/(\frac{\Cdec}{2} +
\Cstab) \, t $) then
 \begin{equation*}
 \left\|f_t - M \right\|_{L^1} \le
 C \left( e^{-\Cstab(t-t_1)}  +
 (t-t_1)^{- \frac{1}{\var} } + e^{- 2 \Cstab (t-t_1) } \right)
 \end{equation*}
and so
 \begin{equation*}
 \left\|f_t - M \right\|_{L^1} \le
 C (t-t_1)^{- \frac{1}{\var} } \le C' t^{- \frac{1}{\var}}.
 \end{equation*}
This holds for $\var$ arbitrarily small, and the theorem is proved.
\end{proof}

%%%%%%%%%%%%%%%%%%%%%%%%%%%%%%%%%%%%%%%%%%%%%%%%%%%%
%%%%%%%%%%%%%%%%%%%%%%%%%%%%%%%%%%%%%%%%%%%%%%%%%%%%
%%%%%%%%%%%%%%%%%%%%%%%%%%%%%%%%%%%%%%%%%%%%%%%%%%%%

\section{Weaker integrability conditions} 
\setcounter{equation}{0}

A natural question is wether the two main results of this paper, the decomposition 
theorem~\ref{theo:dec} and the Theorem~\ref{theo:cvg:eq} 
of convergence to equilibrium, extend to solutions with weaker integrability 
conditions. A first step could be $L^1 _2 \cap L^p$ 
with $1 < p < 2$. A physically relevant assumption would be 
$L^1 _2 \cap L \log L$. 
But since Mischler and Wennberg~\cite{MiscWenn:exun:99} have proven the existence 
and unicity under the sole $L^1 _2$ assumption, the optimal assumption would be only 
$f_0 \in L^1 _2$ (i.e no entropy condition). 

It turns out that in the particular case of hard sphere collision kernel we can extend our  
results to general $L^1 _2$ data, using results of~\cite{MiscWenn:exun:99} and~\cite{Abra:99}. 
A careful study of the iterated gain term $Q^+ ( Q^+ (g,f) , h )$ 
is done in~\cite{MiscWenn:exun:99} in order to prove non-concentration of the solution. This non-concentration 
is used to obtain the weak compactness by Dunford-Pettis Theorem and prove the existence of solution with no 
entropy condition. 
This study is refined in~\cite{Abra:99}, where this iterated gain term is estimated in Lebesgue spaces. 
Therefore Abrahamsson is able to prove~\cite[Lemma~2.1]{Abra:99} 
 \[ \forall \ 1 \le p < 3, \ \ \| Q^+ ( Q^+ (g,f) , h ) \|_{L^p} \le C \, \|f\|_{L^1 _2} 
                                      \, \|g\|_{L^1 _2} \, \|h\|_{L^1 _2} \]
with explicit constant. He deduces a decomposition 
theorem~\cite[Proposition~2.1]{Abra:99} from which we can extract

\begin{lemma}[Abrahamsson's decomposition] \label{lem:abra}
Let $B(v-v_*,\sigma) = |v-v_*|$, and let $f_0\in L^1_2$ be a 
nonnegative initial datum with finite kinetic energy. Let $f$ be the 
unique solution (with non-increasing energy) of the Boltzmann equation 
with collision kernel $B$ and initial datum $f$.
Let $q$ be arbitrarily large and 
$\tau$ arbitrarily small. Then $f$ can be decomposed as $f= f^S + f^R$  
where $f^S \in L^\infty ([\tau,+\infty); L^2 _q \cap L^1 _2)$ and for all $k \ge 0$, 
there is $\lambda = \lambda(k) >0$ such that 
$ \| f^R \|_{L^1 _k} = O(e^{-\lambda t})$. All the constants in this lemma 
can be computed explicitely in term of the mass and energy of $f_0$. 
\end{lemma}

We explain how to connect this result to our method in order to get optimal assumptions on the initial 
data in the hard sphere case, and then we make some remarks on possible extensions for general 
hard potentials with cut-off. 

Thus for hard spheres we have the

 \begin{theorem} \label{theo:SD}
 Let $B(v-v_*, \sigma) = |v-v_*|$ and $0 \le f_0 \in L^1 _2$. 
 Let $f$ be the unique energy-preserving
 solution of the Botzmann equation with initial datum $f_0$, and let
 $s\geq 0$, $q\geq 0$ be arbitrarily large.
 Let $\tau>0$ be arbitrarily small. Then, for all $t\geq \tau$,
 $f$ can be written $f^S+f^R$, where $f^S$ is nonnegative,
 and
  \begin{equation*}
  \left\{
   \begin{array}{l}
   \displaystyle \sup _{t \ge \tau} \left\|f_t ^S \right\|_{H^s_q \cap L^1 _2} < +\infty \\ \\
   \displaystyle \forall t \ge \tau,\> \forall k>0,\> \exists
   \lambda=\lambda(k)>0; \ \ \left\|f_t ^R \right\|_{L^1 _k}
   = O\left(e^{-\lambda t}\right) .
   \end{array}
  \right.
  \end{equation*}
 Moreover the conclusion of Theorem~\ref{theo:cv:2001} holds true:
  \begin{equation*}
  \left\|f_t - M \right\|_{L^1} = O ( t^{-\infty}).
  \end{equation*}
 All the constants in this theorem can be computed in terms of $\tau$
 and the mass and energy of $f_0$.
 \end{theorem}

\begin{proof}[Proof of Theorem~\ref{theo:SD}] 
First let us prove the decomposition part of the theorem. One follows the same 
strategy of tree decomposition as in Theorem~\ref{theo:dec}. It is enough to take the decomposition of 
Lemma~\ref{lem:abra} at the first step of the tree: $f_1$ takes the smooth part of decomposition of 
Lemma~\ref{lem:abra} as initial data at time $t_0$. Then one has to adjust the constants in the proof: 
take $n$, the number of steps, such that $(n+1) \alpha' \ge k$ (one step more) and take 
 \[ \mu < \frac{\Cstab}{\Cstab + K''} \] 
where $K'' = \min \{ K', \lambda \}$ ($\lambda$ is the rate of exponential decrease in the decomposition of 
Lemma~\ref{lem:abra}). 
The rest of the proof is identical to the one of Theorem~\ref{theo:dec}. 
Then with the decomposition result in hand, one can prove the ``almost exponential'' convergence 
to equilibrium eactly the same way as in Theorem~\ref{theo:cvg:eq}.  
\end{proof}

\Remarks 1. Note that except for the physically relevant case of hard spheres,
the cut-off assumption is unphysical for general hard potentials 
interactions. Besides, non cut-off collision operators are known 
to have a regularizing effect (see for instance 
Alexandre, Desvillettes, Villani and Wennberg~\cite{ADVW:00}). 
The optimality of the integrability condition is thus less important 
for general hard potentials interactions than it is for hard spheres. 
\smallskip

2. For general hard potentials with cut-off (with $0 \le \gamma \le 1$), 
the result of Abrahamsson on the iterated gain term becomes 
   \[ \forall \ 1 \le p < 3, \ \ \| Q^+ ( Q^+ (g,f) , h ) \|_{L^p} \le \
                                    C_{p, \gamma} \left(  \|f\|_{L^1 _2}, \|g\|_{L^1 _2}, \|h\|_{L^1 _2}, 
                                    \|f\|_{L^q}, \|g\|_{L^r}, \|b\|_{L^\infty} \right) \]
for any $1/q + 1/r < (5+\gamma)/2$. It is likely that an improvement of this result in order to 
allow $q=r=1$ in this estimate would allow to extend Lemma~\ref{lem:abra} and thus Theorem~\ref{theo:SD} to 
general hard potentials with cut-off. However it seems that this question
leads to serious technical difficulties. 
\smallskip

3. Nevertheless a possible strategy to extend Theorem~\ref{theo:dec} to initial data in $L^1 _2 \cap L^p$ with 
any $p>1$ could be the following. In the same spirit as the tree decomposition in Theorem~\ref{theo:dec}, 
one iterates the Duhamel formula, but now to increase the Lebesgue integrability at each step (using 
Theorem~\ref{theo:reg:Q^+:full} for $s=0$, translated into a gain of integrability thanks to the Sobolev 
injections coupled with some interpolation). As soon as the $L^2$ integrability is reached, one can start 
the decomposition tree of Theorem~\ref{theo:dec} in order to increase regularity, connecting the two decompositions 
in the same spirit as in the proof of Theorem~\ref{theo:SD}.  
\medskip

%%%%%%%%%%%%%%%%%%%%%%%%%%%%%%%%%%%%%%%%%%%%%%%%%%%%
%%%%%%%%%%%%%%%%%%%%%%%%%%%%%%%%%%%%%%%%%%%%%%%%%%%%
%%%%%%%%%%%%%%%%%%%%%%%%%%%%%%%%%%%%%%%%%%%%%%%%%%%%

%%%%%%%%%%%%%%%%%%%%%%%%%%%%%%%%%%%%%%%%%%%%%%%%%%%%
%%%%%%%%%%%%%%%%%%%%%%%%%%%%%%%%%%%%%%%%%%%%%%%%%%%%
%%%%%%%%%%%%%%%%%%%%%%%%%%%%%%%%%%%%%%%%%%%%%%%%%%%%

\bigskip
\appendix
\def\theequation {A.\arabic{equation}}
\def\thetheorem {{A.\arabic{theorem}}}
\def\thesection{}
\section{Appendix: Some facts from
interpolation theory and harmonic analysis}

The goal of this appendix is to recall some classical results
about linear interpolation theory and also to give the proof of
some elementary results used here, in order to make this paper
almost self-contained.

%%%%%%%%%%%% Appendice : section 1
\subsection*{Convolution inequalities in weighted spaces}
\setcounter{equation}{0}
 \begin{proposition}\label{theo:ineqconv}
 Let $\eta \in \R$, then
  \begin{equation*}
  \left\|f * g\right\|_{L^r _\eta} \le
  \left\|f\right\|_{L^p _{|\eta|}}
  \left\|g\right\|_{L^q _\eta}
  \end{equation*}
 for all $p,q,r \ge 1$ such that $\frac{1}{r} +1 = \frac{1}{p} +\frac{1}{q}$.
 \end{proposition}

The proof of this proposition is exactly similar to the standard proof
of the usual Young inequality.

%%%%%%%%%%%% Appendice : section 2
\subsection*{Riesz-Thorin interpolation}
\setcounter{equation}{0}

 \begin{proposition}\label{theo:RT}
 Let $\theta \in [0,1]$, $p_1,p_2,p \in [1,+\infty]$ such that
 $1/p=\theta/p_1+(1-\theta)/p_2$,  $k_1,k_2,k \in \R$ such that
 $k=\theta k_1 + (1-\theta) k_2$, $q_1, q_2, q \in [1,+\infty]$ such
 that $1/q=\theta/q_1+(1-\theta)/q_2$,  $l_1,l_2,l \in \R$
 such that $l=\theta l_1 + (1-\theta) l_2$, and let $T$ be a continuous operator
 from $L^{p_1} _{k_1}$ into $L^{q_1} _{l_1}$ and from $L^{p_2} _{k_2}$
 into $L^{q_2} _{l_2}$.Then its restrictions to $C^\infty _0$ functions
 extends to a continuous operator from $L^{p} _{k}$ into $L^{q} _{l}$ with the
 following bound on its norm
  \begin{equation*}
  \left\|T\right\|_{L^p _k \rightarrow L^q _l} \le
  \left\|T\right\|_{L^{p_1} _{k_1} \rightarrow L^{q_1} _{l_1}}^\theta
  \left\|T\right\|_{L^{p_2} _{k_2} \rightarrow L^{q_2} _{l_2}}^{1-\theta}.
  \end{equation*}
 \end{proposition}

 \begin{corollary}\label{theo:RT:H^s}
 Let $\theta \in [0,1]$, $s_1,s_2,s \in \R$ such that $s=\theta s_1 + (1-\theta) s_2$,
 and $t_1, t_2, t \in \R$ such that $t=\theta t_1 + (1-\theta) t_2$.
 If $T$ is a continuous operator from $H^{s_1}$ into $H^{t_1}$ and from $H^{s_2}$
 into $H^{t_2}$, then its restriction to $C^\infty _0$ functions extends
 to a continuous operator from $H^s$ into $H^t$ with the following bound on its norm
  \begin{equation*}
  \left\|T\right\|_{H^s \rightarrow H^t} \le
  \left\|T\right\|_{H^{s_1} \rightarrow H^{t_1}}^\theta
  \left\|T\right\|_{H^{s_2} \rightarrow H^{t_2}}^{1-\theta}.
  \end{equation*}
 \end{corollary}

This corollary is still true when one adds a weight on the
space variable:
 \begin{equation*}
 \left\|T\right\|_{H^s _k \rightarrow H^t _{k'}} \le
 \left\|T\right\|_{H^{s_1} _k \rightarrow H^{t_1} _{k'}}^\theta
 \left\|T\right\|_{H^{s_2} _k \rightarrow H^{t_2} _{k'}}^{1-\theta}.
 \end{equation*}
In fact the abstract method of interpolation leads to the stronger
result
 \begin{equation*}
 \left\|T\right\|_{H^s_k \rightarrow H^t _{k'}} \le
 \left\|T\right\|_{H^{s_1} _{k_1} \rightarrow H^{t_1} _{k' _1}}^\theta
 \left\|T\right\|_{H^{s_2} _{k_2} \rightarrow H^{t_2} _{k' _2}}^{1-\theta}
 \end{equation*}
where the weight indexes satisfy $k = \theta k_1 + (1-\theta)k_2$
and $k' = \theta k' _1 + (1-\theta)k' _2$.
As a consequence
one could prove a strong version of the Young inequality in the
case of the weighted Sobolev spaces. One can indeed make the index
of weight and regularity vary together. Namely
 \begin{equation*}
 \left\|f\right\|_{H^s _k} \le
 \left\|f\right\|_{H^{s_1} _{k_1}}^\theta
 \left\|f\right\|_{H^{s_2} _{k_2}}^{1-\theta}
 \end{equation*}
where $s = \theta s_1 + (1-\theta)s_2$ and $k = \theta k_1 +
(1-\theta)k_2$. Let us emphasize the consequence of this
inequality that we use in this paper: as soon as $f$ belongs to
$H^s$ and has finite $L^1$ moments of order large enough, one can
deduce bounds on $H^{s'} _k$ norm for $s' < s$.

%%%%%%%%%%%%%%% Appendice : section 3
\subsection*{Regularity of a sum $H^s _k + H^{s+\beta} _k$}
\setcounter{equation}{0}

 \begin{theorem}\label{theo:interp:sum:gen}
 Let $h \in H^s _\eta$ ($s \ge 0$, $\eta \in \R$) such that for all $\var$ small enough,
  \begin{equation*}
  h = h^\var _1 + h^\var _2
  \end{equation*}
 where the two parts $h^\var _1$ and $h^\var _2$ satisfy the
 following estimates: there exist $k_1 \ge 0$ and $k_2 >0$ such that
 $\left\|h^\var _1\right\|_{H^{s+\beta} _\eta} \le C_1 \var ^{-k_1}$ et
 $\left\|h^\var _2\right\|_{H^s _\eta} \le C_2 \var ^{k_2}$
($\beta >0$).
 Then
  \begin{equation*}
  h \in H^{s+\alpha} _\eta,\ \  \forall \alpha <  \frac{\beta k_2}{k_1 + k_2}.
  \end{equation*}
\end{theorem}

\Remarks 1.Our estimate on the norm $H^{s + \alpha} _\eta$ blows up like
 $\frac{\cst}{\frac{\beta k_2}{k_1 + k_2} - \alpha}$ as
 $\alpha \to \frac{\beta k_2}{k_1 + k_2}$.
\smallskip

2. In fact the proof shows that $h\langle \cdot\rangle^\eta$ belongs to the Besov
space $B_\beta ^{\infty, \infty}\subset \cap _{\alpha < \beta} H_\alpha$.
\medskip

\begin{proof}[Proof of theorem~\ref{theo:interp:sum:gen}]
Let us take $\alpha < \frac{\beta k_2}{k_1 + k_2}$. Without loss of generality we
treat the case $s = \eta=0$ (the general case can be reduced to this one). We first
prove an upper bound on an annulus. Let $0 < A < B$, and
 \begin{eqnarray*}
 \int_{A \le |\xi| \le B} |\widehat{h}(\xi)|^2 \langle \xi \rangle^{2\alpha} d\xi
 &\le&
 2 \int_{A \le |\xi| \le B} \left( |\widehat{h^\var _1}(\xi)|^2+
 |\widehat{h^\var _2}(\xi)|^2 \right) \langle \xi \rangle^{2\alpha} d\xi \\ 
 &\le& 2 \left( C_1 \var ^{-k_1} \langle A \rangle^{2(\alpha - \beta)} +
 C_2 \langle B \rangle ^{2\alpha} \var ^{k_2} \right) \\ 
 &\le& 2 \left( C_1 \var ^{-k_1}A^{2(\alpha - \beta)} +
 2 C_2 B^{2 \alpha} \var ^{k_2} \right) \\ 
 &\le& \max(2 C_1, 4 C_2) \left( \var ^{-k_1}A^{2(\alpha - \beta)} +
 B^{2 \alpha} \var ^{k_2} \right)
 \end{eqnarray*}
As this inequality holds for all $\var$, one can choose it in order that the two
right-members be equal in the preceding inequality. The computation leads to
 \begin{equation*}
 \int_{A \le |\xi| \le B} |\widehat{h}(\xi)|^2 \langle \xi \rangle^\beta d\xi \le
 2 \max(2 C_1, 4 C_2) B^{\frac{2\alpha k_1}{k_1+k_2}} A^{\frac{2(\alpha - \beta) k_2}{k_1+k_2}}
 \end{equation*}
Let $C_3=2 \max(2 C_1, 4 C_2)$ and let us sum the inequalities
on a family of concentric dyadic annuli:
 \begin{eqnarray*}
 \left\|h\right\|_{H^\beta} &\le&
 \int_{0 \le |\xi| \le 1} |\widehat{h}(\xi)|^2 \langle \xi \rangle^\beta
 + C_3 \sum_{n=0}^{+\infty}2^{\frac{2\alpha (n+1) k_1}{k_1+k_2}}
 2^{\frac{2(\alpha - \beta)n k_2}{k_1+k_2}}\\ 
 &\le& 2 \left\|h\right\|_{L^2} +
 C_3 \sum_{n=0}^{+\infty}2^{\frac{2\alpha (n+1) k_1}{k_1+k_2}}
 2^{\frac{2(\alpha - \beta)n k_2}{k_1+k_2}}\\ 
 &\le& 2 \left\|h\right\|_{L^2} + C_3 4^{\frac{\alpha k_1}{k_1+k_2}}
 \sum_{n=0}^{+\infty}4^{n(\alpha - \frac{\beta k_2}{k_1+k_2})}.
 \end{eqnarray*}
Thanks to the assumption on $\alpha$ the right member is summable and thus
$h \in H^\alpha$ with the following bound on the norm
 \begin{eqnarray*}
 \left\|h\right\|_{H^\alpha} &\le&
 2 \left\|h\right\|_{L^2} +
 C_3 4^{\frac{\alpha k_1}{k_1+k_2}} \sum_{n=0}^{+\infty}
 4^{n(\alpha - \frac{\beta k_2}{k_1+k_2}}) \\ 
 &\le& 2 \left\|h\right\|_{L^2} +
 C_3 \frac{4^{\frac{\alpha k_1}{k_1+k_2}}}{1-4^{\alpha - \frac{\beta k_2}{k_1+k_2}}}.
 \end{eqnarray*}
\end{proof}

%%%%%%%%%%%%%%% Appendice : section 4
\subsection*{A simple estimate on pseudo-differential operators}
\setcounter{equation}{0}

We conclude this appendix with a simple result needed for the
proof of the regularity property of $Q^+$. This can be linked with
more general pseudo-differential estimates, but will be proved by
elementary means. The space $H^s _k$ is
not an algebra in general (it is an algebra
thanks to the Sobolev imbeddings as soon as $2s >N$),
but one can prove a bound on the norm $H^s _k$ of a
product of functions if one of the two functions has regularity
greater than $s$:
 \begin{equation*}
 \left\|fg\right\|_{H^s _k} \le \cst(N,\var) \, \left\|f\right\|_{H^S _{k_1}}
 \, \left\|g\right\|_{H^s _{k_2}}
 \end{equation*}
where $k_1 + k_2=k$, and $S = s +  N/2  + \var$ with $\var > 0$.

Now we follow the same idea but assuming that one of the two
functions depends also on the Fourier variable.
 \begin{lemma}\label{lem:OPD}
 Let $\psi (x,\xi)$ be a real-valued $C^\infty$ function
 on $\R^N\times\R^N$, compactly supported in $x$, uniformly in $\xi$.
 Let $g : \R^N \longrightarrow \R$ be a function in the Schwartz space
 $\cal{S}(\R^N)$, and let $s \in \R$.
 Let us define
  \begin{equation*}
  I = \int_{\R^N} \langle\xi\rangle^{2s}
  \Bigl | \mathcal{F} (g(\cdot)\psi(\cdot,\xi)) \Bigr |^2 d\xi.
  \end{equation*}
 Then for all $\var >0$
 there exists a constant $\cst(N,\var)$ such that
  \begin{equation*}
  I \le \cst (N,\var) \left\|\psi\right\|_{L^\infty _\xi (H^S _x)} ^2
  \left\|g\right\|_{H^s} ^2,
  \end{equation*}
 with $S = s + N/2 + \var$.
 \end{lemma}

\begin{proof}[Proof of Lemma~\ref{lem:OPD}]
We have
 \begin{equation*}
 g(x) = \frac{1}{(2\pi)^{N/2}}
 \int_{\R^N} e^{ix \cdot \tau} \hat{g}(\tau) d\tau
 \end{equation*}
hence
 \begin{eqnarray*} 
 I &=& \frac{1}{(2\pi)^N}
 \int_{\R^N} \langle\xi\rangle^{2s}
 \left| \int_{\R^N} \int_{\R^N} e^{-ix \cdot (\xi - \tau)}
 \hat{g}(\tau) \psi(x,\xi) dx d\tau \right|^2 d\xi \\ \nonumber
 &=& \frac{1}{(2\pi)^N}\int_{\R^N} \langle\xi\rangle^{2s}
 \left| \int_{\R^N} \hat{g}(\tau)
 \left[\int_{\R^N} e^{-ix \cdot (\xi - \tau)}
 \psi(x,\xi) dx \right] d\tau \right|^2 d\xi \\ \nonumber
 &=& \frac{1}{(2\pi)^N}\int_{\R^N} \langle\xi\rangle^{2s}
 \left| \int_{\R^N} \hat{g}(\xi - \tau)
 \left[\int_{\R^N} e^{-ix \cdot \tau}
 \psi(x,\xi) dx \right] d\tau \right|^2 d\xi \\
 &=& \frac{1}{(2\pi)^N}\int_{\R^N} \langle\xi\rangle^{2s}
 \left| \int_{\R^N} \hat{g}(\xi - \tau)
 \left[\mathcal{F} _x
 (\psi(\cdot,\xi)) (\tau) \right] d\tau \right|^2 d\xi
 \end{eqnarray*}
and thus
 \begin{eqnarray} \nonumber
 (2\pi)^N I &\le& \int_{\R^N} \langle\xi\rangle^{2s}
 \int_{\R^N} |\hat{g}|^2(\xi - \tau)
 \langle\tau\rangle^{-2S} d\tau  \\ \nonumber
 && \qquad
 \int_{\R^N}
 \bigl | \mathcal{F} _x (\psi(\cdot,\xi)) \bigr |^2 (\tau ')
 \langle\tau'\rangle^{2S} d\tau' d\xi \\ \nonumber
 &\le& \int_{\R^N} \langle\xi\rangle^{2s}
 \int_{\R^N} | \hat{g}|^2(\xi - \tau) \langle\tau\rangle^{-2S} d\tau
 \left\|\psi(\cdot,\xi)\right\|_{H^S _x}^2
 d\xi \\ \nonumber
 &\le& \left\|\psi\right\|_{L^\infty _\xi (H^S _x)} ^2
 \int_{\R^N} \langle\tau\rangle^{-2S}
 \int_{\R^N} \langle\xi\rangle^{2s} |\hat{g}|^2(\xi - \tau) d\xi\, d\tau
 \\ \nonumber
 &\le& \left\|\psi\right\|_{L^\infty _\xi (H^S _x)} ^2
 \int_{\R^N}  \langle\tau\rangle^{-2S}
 \int_{\R^N} \langle\xi+\tau\rangle^{2s} |\hat{g}(\xi)|^2 d\xi\, d\tau  \\ \nonumber
 &\le& \left\|\psi\right\|_{L^\infty _\xi (H^S _x)} ^2
 \, 2^s \int_{\R^N} \langle\xi\rangle^{2s} \hat{g}(\xi)^2 d\xi
 \int_{\R^N} \langle\tau\rangle^{-2S} \langle\tau\rangle^{2s} d\tau \\ \nonumber
 &\le& \left\|\psi\right\|_{L^\infty _\xi (H^S _x)} ^2
 \left\|g\right\|_{H^s} ^2 \, 2^s \int_{\R^N}
 \langle\tau\rangle^{- N - 2 \var} d\tau
 \end{eqnarray}
which concludes the proof.
\end{proof}

\bigskip

\noindent 
{\bf{Acknowledgment}}: The authors thank the referee for useful comments. 
Support by the European network HYKE, funded by the EC as
contract HPRN-CT-2002-00282, is acknowledged. 
% HYKE
\bigskip

\bibliographystyle{acm}

\bibliography{./SHEB}

\begin{thebibliography}{10}

\bibitem{Abra:99}
{\sc Abrahamsson, F.}
\newblock Strong {$L\sp 1$} convergence to equilibrium without entropy
  conditions for the {B}oltzmann equation.
\newblock {\em Comm. Partial Differential Equations 24}, 7-8 (1999),
  1501--1535.

\bibitem{ADVW:00}
{\sc Alexandre, R., Desvillettes, L., Villani, C., and Wennberg, B.}
\newblock Entropy dissipation and long-range interactions.
\newblock {\em Arch. Ration. Mech. Anal. 152}, 4 (2000), 327--355.

\bibitem{Arke:I+II:72}
{\sc Arkeryd, L.}
\newblock On the {B}oltzmann equation.
\newblock {\em Arch. Rational Mech. Anal. 45\/} (1972), 1--34.

\bibitem{Arke:infty:81}
{\sc Arkeryd, L.}
\newblock Intermolecular forces of infinite range and the {B}oltzmann equation.
\newblock {\em Arch. Rational Mech. Anal. 77\/} (1981), 11--21.

\bibitem{Arke:Linfty:83}
{\sc Arkeryd, L.}
\newblock {$L^\infty$} estimates for the space-homogeneous {B}oltzmann
  equation.
\newblock {\em J. Statist. Phys. 31}, 2 (1983), 347--361.

\bibitem{Arke:stab:88}
{\sc Arkeryd, L.}
\newblock Stability in {$L\sp 1$} for the spatially homogeneous {B}oltzmann
  equation.
\newblock {\em Arch. Rational Mech. Anal. 103}, 2 (1988), 151--167.

\bibitem{Boby:maxw:88}
{\sc Bobyl{\"e}v, A.~V.}
\newblock The theory of the nonlinear spatially uniform {B}oltzmann equation
  for {M}axwell molecules.
\newblock In {\em Mathematical physics reviews, Vol.\ 7}. Harwood Academic
  Publ., Chur, 1988, pp.~111--233.

\bibitem{BoucDesv:rgain:98}
{\sc Bouchut, F., and Desvillettes, L.}
\newblock A proof of the smoothing properties of the positive part of
  {B}oltzmann's kernel.
\newblock {\em Rev. Mat. Iberoamericana 14}, 1 (1998), 47--61.

\bibitem{Carl:fond:32}
{\sc Carleman, T.}
\newblock Sur la th\'eorie de l'equation int\'egrodiff\'erentielle de
  {B}oltzmann.
\newblock {\em Acta Math. 60\/} (1932), 369--424.

\bibitem{Carl:fond:57}
{\sc Carleman, T.}
\newblock {\em Probl\`emes Math\'ematiques dans la Th\'eorie Cin\'etique des
  Gaz}.
\newblock Almqvist \& Wiksell, 1957.

\bibitem{Desv:momt:93}
{\sc Desvillettes, L.}
\newblock Some applications of the method of moments for the homogeneous
  {B}oltzmann and {K}ac equations.
\newblock {\em Arch. Rational Mech. Anal. 123}, 4 (1993), 387--404.

\bibitem{DesvVill:LChomII:2000}
{\sc Desvillettes, L., and Villani, C.}
\newblock On the spatially homogeneous {L}andau equation for hard potentials.
  {I}{I}. ${H}$-theorem and applications.
\newblock {\em Comm. Partial Differential Equations 25}, 1-2 (2000), 261--298.

\bibitem{Goud:nocu:97}
{\sc Goudon, T.}
\newblock On {B}oltzmann equations and {F}okker-{P}lanck asymptotics: influence
  of grazing collisions.
\newblock {\em J. Statist. Phys. 89}, 3-4 (1997), 751--776.

\bibitem{Grad:58}
{\sc Grad, H.}
\newblock Principles of the kinetic theory of gases.
\newblock In {\em {F}l{\"{u}}gge's {H}andbuch des {P}hysik}, vol.~XII.
  Springer-Verlag, 1958, pp.~205--294.

\bibitem{Gust:L^p:86}
{\sc Gustafsson, T.}
\newblock ${L}\sp p$-estimates for the nonlinear spatially homogeneous
  {B}oltzmann equation.
\newblock {\em Arch. Rational Mech. Anal. 92}, 1 (1986), 23--57.

\bibitem{Gust:L^p:88}
{\sc Gustafsson, T.}
\newblock Global ${L}\sp p$-properties for the spatially homogeneous
  {B}oltzmann equation.
\newblock {\em Arch. Rational Mech. Anal. 103}, 1 (1988), 1--38.

\bibitem{Lion:rgainI+II:94}
{\sc Lions, P.-L.}
\newblock Compactness in {B}oltzmann's equation via {F}ourier integral
  operators and applications. {I}, {I}{I}.
\newblock {\em J. Math. Kyoto Univ. 34}, 2 (1994), 391--427, 429--461.

\bibitem{Lion:rgainIII:94}
{\sc Lions, P.-L.}
\newblock Compactness in {B}oltzmann's equation via {F}ourier integral
  operators and applications. {I}{I}{I}.
\newblock {\em J. Math. Kyoto Univ. 34}, 3 (1994), 539--584.

\bibitem{Lu:rgain:98}
{\sc Lu, X.}
\newblock A direct method for the regularity of the gain term in the
  {B}oltzmann equation.
\newblock {\em J. Math. Anal. Appl. 228}, 2 (1998), 409--435.

\bibitem{MiscWenn:exun:99}
{\sc Mischler, S., and Wennberg, B.}
\newblock On the spatially homogeneous {B}oltzmann equation.
\newblock {\em Ann. Inst. H. Poincar\'e Anal. Non Lin\'eaire 16}, 4 (1999),
  467--501.

\bibitem{Povz:ineq:65}
{\sc Povzner, A.~J.}
\newblock The {B}oltzmann equation in the kinetic theory of gases.
\newblock {\em Amer. Math. Soc. Transl. 47}, Ser. 2 (1965), 193--214.

\bibitem{PulvWenn:binf:97}
{\sc Pulvirenti, A., and Wennberg, B.}
\newblock A {M}axwellian lower bound for solutions to the {B}oltzmann equation.
\newblock {\em Comm. Math. Phys. 183}, 1 (1997), 145--160.

\bibitem{SoggStei:radon:85}
{\sc Sogge, C.~D., and Stein, E.~M.}
\newblock Averages of functions over hypersurfaces in {${\bf {R}}\sp n$}.
\newblock {\em Invent. Math. 82}, 3 (1985), 543--556.

\bibitem{SoggStei:radon:86}
{\sc Sogge, C.~D., and Stein, E.~M.}
\newblock Averages over hypersurfaces. {I}{I}.
\newblock {\em Invent. Math. 86}, 2 (1986), 233--242.

\bibitem{SoggStei:radon:90}
{\sc Sogge, C.~D., and Stein, E.~M.}
\newblock Averages over hypersurfaces. {S}moothness of generalized {R}adon
  transforms.
\newblock {\em J. Analyse Math. 54\/} (1990), 165--188.

\bibitem{ToscVill:cveq:2000}
{\sc Toscani, G., and Villani, C.}
\newblock On the trend to equilibrium for some dissipative systems with slowly
  increasing a priori bounds.
\newblock {\em J. Statist. Phys. 98}, 5-6 (2000), 1279--1309.

\bibitem{Vill:nocu:98}
{\sc Villani, C.}
\newblock On a new class of weak solutions to the spatially homogeneous
  {B}oltzmann and {L}andau equations.
\newblock {\em Arch. Rational Mech. Anal. 143}, 3 (1998), 273--307.

\bibitem{Vill:cveq:TA}
{\sc Villani, C.}
\newblock Cercignani's conjecture is sometimes true, and always almost true.
\newblock Preprint, 2002.

\bibitem{Wenn:stab:93}
{\sc Wennberg, B.}
\newblock Stability and exponential convergence in {$L\sp p$} for the spatially
  homogeneous {B}oltzmann equation.
\newblock {\em Nonlinear Anal. 20}, 8 (1993), 935--964.

\bibitem{Wenn:momt:94}
{\sc Wennberg, B.}
\newblock On moments and uniqueness for solutions to the space homogeneous
  {B}oltzmann equation.
\newblock {\em Transport Theory Statist. Phys. 23}, 4 (1994), 533--539.

\bibitem{Wenn:rado:94}
{\sc Wennberg, B.}
\newblock Regularity in the {B}oltzmann equation and the {R}adon transform.
\newblock {\em Comm. Partial Differential Equations 19}, 11-12 (1994),
  2057--2074.

\bibitem{Wenn:momt:97}
{\sc Wennberg, B.}
\newblock Entropy dissipation and moment production for the {B}oltzmann
  equation.
\newblock {\em J. Statist. Phys. 86}, 5-6 (1997), 1053--1066.

\bibitem{Wenn:nonu:99}
{\sc Wennberg, B.}
\newblock An example of nonuniqueness for solutions to the homogeneous
  {B}oltzmann equation.
\newblock {\em J. Statist. Phys. 95}, 1-2 (1999), 469--477.

\end{thebibliography}

\begin{flushleft} \signcm \end{flushleft}
\vspace*{-55mm} \begin{flushright} \signcv \end{flushright}

\end{document}